%% file: sigma24-071.tex
\pdfoutput=1
\RequirePackage{ifpdf}
\ifpdf 
\documentclass[pdftex]{sigma}
\else
\documentclass{sigma}
\fi

\numberwithin{equation}{section}

\newtheorem{Theorem}{Theorem}[section]
\newtheorem{Corollary}[Theorem]{Corollary}
\newtheorem{Lemma}[Theorem]{Lemma}

 { \theoremstyle{definition}
\newtheorem{Definition}[Theorem]{Definition}

\newtheorem{Example}[Theorem]{Example}
\newtheorem{Remark}[Theorem]{Remark} }

\usepackage{mathtools}
\newcommand{\Span}[1]{\mathrm{span}\{#1\}}
\newcommand{\Dim}[1]{\dim(#1)}

\newcommand{\D}{\mathrm{d}}
\newcommand{\Lie}{\mathrm{L}}

\newcommand{\Rank}[1]{\mathrm{rank}(#1)}

\newcommand{\XxUlu}{\mathcal{X}\times\mathcal{U}_{[0,l_u]}}
\usepackage{lineno}

\begin{document}

\allowdisplaybreaks

\newcommand{\arXivNumber}{2110.12995}

\renewcommand{\PaperNumber}{071}

\FirstPageHeading

\ShortArticleName{Tracking Control for $(x,u)$-Flat Systems by Quasi-Static Feedback of Classical States}

\ArticleName{Tracking Control for $\boldsymbol{(x,u)}$-Flat Systems\\ by Quasi-Static Feedback of Classical States}
\Author{Conrad GST\"OTTNER, Bernd KOLAR and Markus SCH\"OBERL}
\AuthorNameForHeading{C.~Gst\"ottner, B.~Kolar and M.~Sch\"oberl}
\Address{Institute of Automatic Control and Control Systems Technology,\\ Johannes Kepler University Linz, Austria}
\Email{\href{mailto:conrad.gstoettner@jku.at}{conrad.gstoettner@jku.at}, \href{mailto:bernd_kolar@ifac-mail.org}{bernd\_kolar@ifac-mail.org}, \href{mailto:markus.schoeberl@jku.at}{markus.schoeberl@jku.at}}

\ArticleDates{Received November 03, 2023, in final form July 22, 2024; Published online July 31, 2024}

\Abstract{It is well known that for flat systems the tracking control problem can be solved by utilizing a linearizing quasi-static feedback of generalized states. If measurements (or estimates) of a so-called generalized Brunovsk\'y state are available, a linear, decoupled and asymptotically stable tracking error dynamics can be achieved. However, from a practical point of view, it is often desirable to achieve the same tracking error dynamics by feedback of a~classical state instead of a~generalized one. This is due to the fact that the components of a~classical state typically correspond to measurable physical quantities, whereas a generalized Brunovsk\'y state often contains higher order time derivatives of the (fictitious) flat output which are not directly accessible by measurements. In this paper, a systematic solution for the tracking control problem based on quasi-static feedback and measurements of classical states only is derived for the subclass of $(x,u)$-flat systems.}


\Classification{53Z30; 93C10; 93C15}

\section{Introduction}\label{sec:intro}
The concept of flatness has been introduced in control theory by Fliess, L\'evine, Martin and Rouchon in the 1990s, see, e.g., \cite{FliessLevineMartinRouchon:1992-2,FliessLevineMartinRouchon:1992,FliessLevineMartinRouchon:1995}. Roughly speaking, a nonlinear control system of the form\looseness=-1%
\begin{align}\label{eq:intro_nlsys}
	\dot{x}&=f(x,u)
\end{align}
with an $n$-dimensional\footnote{By abuse of notation, we write, e.g., $\Dim{x}=n$, although it is actually the space containing $x$, which is of dimension $n$. The same abuse of notation is also used for other quantities throughout the paper.} state $x$ and an $m$-dimensional input $u$ is flat, if there exist $m$ differentially independent functions $y^j=\varphi^j\bigl(x,u,\dot u,\dots,u^{(q)}\bigr)$, such that $x$ and $u$ can locally be expressed by~$y$ and its time derivatives, i.e.,
\begin{align}\label{eq:intro_parameterization}
	x=F_x\bigl(y,\dot y,\dots,y^{(r-1)}\bigr),\qquad u=F_u\bigl(y,\dot y,\dots,y^{(r)}\bigr).
\end{align}
Parameterization \eqref{eq:intro_parameterization} of the system variables by the flat output allows for an elegant and systematic solution of both feedforward and feedback problems, which is the reason for the ongoing popularity of flat systems, see, e.g., \cite{DelaleauRudolph:1998,FliessLevineMartinRouchon:1995}. The computation of flat outputs, however, is known to be a difficult problem. Recent research in this field can be found, e.g., in \cite{GstottnerKolarSchoberl:2022,GstottnerKolarSchoberl:2021-3,NicolauRespondek:2017,NicolauRespondek:2019}.

A typical intermediate step in the design of a flatness-based control is an exact linearization of system \eqref{eq:intro_nlsys} by a suitable feedback. Exact feedback linearization should not be confused with exact feedforward linearization, which has been introduced in \cite{HagenmeyerDelaleau:2003}, see also \cite{HagenmeyerDelaleau:2003-2}. In this paper, we solely consider exact feedback linearization. The objective is that the closed-loop system possesses a linear input-output behavior in the form of integrator chains between a new closed-loop input $v$ and the considered flat output $y$. A standard approach is an exact linearization by an endogenous dynamic feedback
\begin{align}\label{eq:intro_dynamic_feedback}
		\dot{z}&=g(x,z,v),\qquad
		u=\alpha(x,z,v),
\end{align}
where the new input $v$ is given by the highest time derivatives $y^{(r)}$ that are present in parameterization \eqref{eq:intro_parameterization}. The feedback \eqref{eq:intro_dynamic_feedback} being endogenous means that $g$ and $\alpha$ are such that the state $z$ of the feedback and the new input $v$ which is introduced by means of this feedback can be expressed as functions of $x$, $u$ and time derivatives of $u$. Thus, the closed-loop system~${\dot{x}\!=f(x,\alpha(x,z,v))}$,~${\dot{z}\!=g(x,z,v)}
$
has a linear input-output behavior~\smash{$y^{j,(r^j)}\!=v^j$}, ${j\!=1,\dots,m}$ between the new input $v$ and the flat output $y$. For example, in \cite{Chetverikov:2007,FliessLevineMartinRouchon:1999,Kolar:2017} it is shown how such a dynamic feedback \eqref{eq:intro_dynamic_feedback} can be constructed systematically from parameterization \eqref{eq:intro_parameterization} (we recall the construction in Section~\ref{sec:exact_lin}).

However, a drawback of this standard approach is that it is a dynamic feedback and that the order of the error dynamics of a subsequently designed tracking control is given by $\dim(x)+\dim(z)$, which is higher than the order $\dim(x)=n$ of the original system \eqref{eq:intro_nlsys}. A well-established alternative, which circumvents these drawbacks, is the exact linearization by a quasi-static feedback of generalized states proposed in \cite{DelaleauRudolph:1995,DelaleauRudolph:1998,RudolphDelaleau:1998}. In contrast to the more classical static or dynamic feedbacks, a quasi-static feedback is a feedback which may also depend on time derivatives of the closed-loop input. The notion of quasi-static feedback and the equivalence of two systems by a quasi-static feedback were originally defined within a differential-algebraic framework, see, e.g., \cite{DelaleauFliess:1992-2,DelaleauFliess:1992,DelaleauRudolph:1995,Rudolph:1995} for a precise definition and further details. A consequence of this definition is the existence of a transformation of the form
	\begin{align}\label{eq:quasistatic_k}
			v^{(k)}=\phi_k(x,u,\dot u,\dots),\qquad
			u^{(k)}=\hat\phi_k(x,v,\dot v,\dots)
	\end{align}
	for $k\geq 0$, which relates the original input $u$ and the new input $v$ as well as their time derivatives.
Note that although $x$ is a classical state for the original system $\dot{x}=f(x,u)$, it is a generalized one for the resulting closed-loop system $\dot{x}=f(x,v,\dot{v},\dots)$.
The fact that one input can be calculated directly from the other by relation \eqref{eq:quasistatic_k}, without the need to solve any differential equation, explains the terminology quasi-static, see also \cite{Rudolph:2021}. In \cite{DelaleauRudolph:1995,DelaleauRudolph:1998}, it has been shown that every flat system \eqref{eq:intro_nlsys} can be exactly linearized by a quasi-static feedback of a so-called generalized Brunovsk\'y state $\tilde{x}_B=\bigl(y,\dot y,\dots,y^{(\kappa-1)}\bigr)$, which consists of suitably chosen time derivatives of the components $y^j$ of the flat output up to the orders $\kappa^j-1$ and meets $\dim(\tilde{x}_B)=\dim(x)$, i.e., $\kappa^1+\dots+\kappa^m=n$. Passing from the original state $x$ to this generalized Brunovsk\'y state~$\tilde{x}_B$ results in a state representation of the form $\dot{\tilde{x}}_B=f(\tilde{x}_B,u,\dot{u},\dots)$, i.e., a generalized state representation. In a next step, the system is exactly linearized by a quasi-static feedback	
\begin{align}\label{eq:intro_quasistatic_feedback}
	u&=\alpha(\tilde{x}_B,v,\dot{v},\dots),
\end{align}
where the new input $v$ is given by the time derivatives $y^{(\kappa)}$, and hence the closed-loop system
possesses the linear input-output behavior $y^{j,(\kappa^j)}=v^j$, $j=1,\dots,m$. Note that here $v$ corresponds to lower-order time derivatives of the flat output than in the dynamic feedback \eqref{eq:intro_dynamic_feedback}. Under the assumption that measurements or estimates of the generalized Brunovsk\'y state $\tilde{x}_B$ are available, it is then straightforward to achieve in a second step a linear and decoupled tracking error dynamics for the components of the flat output with arbitrarily placed eigenvalues. The practical usefulness and good performance which can be achieved with a quasi-static tracking control has been demonstrated, e.g., in \cite{KugiKiefer:2005}.
However, feedback \eqref{eq:intro_quasistatic_feedback} is not a feedback of the original state~$x$ of \eqref{eq:intro_nlsys}, and the same applies to the additional feedback which achieves the desired tracking error dynamics. Since except for static feedback linearizable systems a generalized Brunovsk\'y state $\tilde{x}_B$ is not equivalent to the original state~$x$ via a state transformation $\tilde{x}_B =\Phi_x(x)$, determining the required time derivatives of the flat output from available measurements is not a~straightforward task.

The aim of the present paper is to propose a method which combines the advantages of both above-mentioned approaches while avoiding their drawbacks. The solution of the considered problem requires two steps. First, a linearizing quasi-static feedback
\begin{align}\label{eq:intro_quasi_static_feedback_x}
		&u=\alpha(x,v,\dot{v},\dots)
\end{align}
of the original, classical state $x$ has to be derived. Second, for the exactly feedback linearized system, an additional feedback has to be constructed which achieves the desired linear, decoupled and asymptotically stable tracking error dynamics and again only requires the classical state $x$ of the open-loop system \eqref{eq:intro_nlsys}. Regarding the first step, there are already results available in the literature. In particular, in \cite{DelaleauPereiradaSilva:1998-1,DelaleauPereiradaSilva:1998-2} conditions for the input-output linearizability by a~quasi-static feedback of classical states have been derived in a~differential-algebraic setting and it can be shown that $(x,u)$-flat systems satisfy these conditions. However, the exact linearization by a~quasi-static feedback of the form \eqref{eq:intro_quasi_static_feedback_x} is only an intermediate step in deriving a tracking control law of the desired form $u=\alpha\bigl(x,y^d(t),\dot{y}^d(t),\dots,y^{d,(r)}(t)\bigr)$, which only depends on $x$ and the reference trajectory $y^d(t)$. Therefore, we propose an alternative derivation of a linearizing quasi-static feedback \eqref{eq:intro_quasi_static_feedback_x} within a differential-geometric framework, which as a byproduct yields a coordinate transformation with a certain triangular structure that is essential for the subsequent tracking control design.

At this point, it is important to emphasize that transforming a feedback of the form \eqref{eq:intro_quasistatic_feedback} into a feedback of the form \eqref{eq:intro_quasi_static_feedback_x} is not a straightforward task and not generally possible.
\begin{Remark}
	Since the generalized Brunovsk\'y state $\tilde{x}_B$ consists of certain time derivatives of components of the flat output, it is of course a function of $x$, $u$, and time derivatives of $u$, i.e., $\tilde{x}_B=\phi(x,u,\dot{u},\dots)$. However, simply substituting this expression into \eqref{eq:intro_quasistatic_feedback} yields in general an implicit relation of the form\footnote{With the exception of $y$ being a linearizing output in the sense of static feedback linearizability, at least one of the components of $\tilde{x}_B=\phi(x,u,\dot{u},\dots)$ explicitly depends on $u$ or time derivatives of $u$.}
	\begin{align}\label{eq:quasistatic_feedback_x_implicit}
			u&=\alpha(\phi(x,u,\dot{u},\dots),v,\dot{v},\dots)
	\end{align}
	in $u$ and its time derivatives. To obtain a quasi-static feedback of the desired form \eqref{eq:intro_quasi_static_feedback_x}, we would have to differentiate \eqref{eq:quasistatic_feedback_x_implicit} in order to eventually obtain a system of equations which can be solved for $u$ as a function of $x$, $v$, and time derivatives of $v$. According to the author's best knowledge, it is not known whether for every flat system \eqref{eq:intro_nlsys} there exists at least one linearizing feedback \eqref{eq:intro_quasistatic_feedback} for which this is indeed possible.
\end{Remark}
The contribution of the paper is threefold: First, we derive geometric conditions under which certain time derivatives of a general flat output
$
		y=\varphi(x,u,\dot{u},\dots)
$
can be introduced as a new input, and provide a systematic procedure for the construction of the corresponding feedback. Second, we show within our differential geometric framework that for $(x,u)$-flat systems, i.e., systems with a flat output
\begin{align}\label{eq:xu_flat_output}
	y&=\varphi(x,u),
\end{align}
which may depend on the input but not on time derivatives of the input, the new input can always be chosen such that the corresponding feedback is a quasi-static feedback of the systems' original, classical state $x$. In contrast to \cite{DelaleauPereiradaSilva:1998-1,DelaleauPereiradaSilva:1998-2}, this alternative approach reveals structural properties which are important regarding tracking control design. Third, we show that on the basis of such an exact feedback linearization, it is possible to achieve a linear, decoupled and asymptotically stable tracking error dynamics by an additional feedback which again requires only measurements of the state $x$. Thus, in contrast to the standard approach for the tracking control design described in \cite{DelaleauRudolph:1998}, the usage of a generalized Brunovsk\'y state can again be avoided. Preliminary results addressing this topic can be found in \cite{KolarRamsSchlacher:2017} and \cite{GstottnerKolarSchoberl:2020}. The approach taken therein roughly speaking consists of deriving a tracking control law based on a~generalized Brunovsk\'y state and -- under certain assumptions -- subsequently transforming this control law such that it depends on the original, classical state $x$ and the reference trajectory only. In the present contribution, in contrast, we completely avoid the usage of a generalized Brunovsk\'y state and show for $(x,u)$-flat systems how to derive a tracking control law of the desired form~${u=\alpha\bigl(x,y^d(t),\dot{y}^d(t),\dots,y^{d,(r)}(t)\bigr)}$, which only depends on $x$ and the reference trajectory $y^d(t)$, directly.

The paper is organized as follows: In Section \ref{sec:preliminaries}, we introduce some notation and preliminaries. In Section \ref{sec:exact_lin}, we discuss the exact feedback linearization of flat systems in a differential-geometric framework and show how to construct a feedback which introduces appropriately chosen time derivatives of a flat output as new input. Furthermore, we prove in this framework that for~$(x,u)$-flat systems the new input can always be chosen in such a way that the required feedback is a~quasi-static feedback of the state $x$. The results of Section \ref{sec:exact_lin} are then illustrated by two examples in Section \ref{sec:examples}. Subsequently, Section \ref{sec:control} deals with the design of a flatness-based tracking control on the basis of the exact feedback linearization derived in Section~\ref{sec:exact_lin}. Finally, in Section~\ref{sec:examples_continued}, the tracking control design is illustrated with the continued examples of Section~\ref{sec:examples}.\looseness=1

\section{Preliminaries}\label{sec:preliminaries}
	In the following, some notation and the utilized differential-geometric framework are introduced.\looseness=1

	\subsection{Notation}
		Let $\mathcal{X}$ be an $n$-dimensional smooth manifold, equipped with local coordinates $x^i$, $i=1,\dots,n$. The tangent bundle and the cotangent bundle of $\mathcal{X}$ are denoted by $(\mathcal{T}(\mathcal{X}),\tau_\mathcal{X},\mathcal{X})$ and $(\mathcal{T}^\ast(\mathcal{X}),\allowbreak\tau^\ast_\mathcal{X},\mathcal{X})$. For these bundles we have the induced local coordinates $\bigl(x^i,\dot{x}^i\bigr)$ and $\bigl(x^i,\dot{x}_i\bigr)$ with respect to the holonomic bases $\big\{\partial_{x^i}\big\}$ and $\big\{\D x^i\big\}$, respectively. We also make use of the Einstein summation convention. A vector field is a section of the tangent bundle, i.e., a map ${w\colon\mathcal{X}\rightarrow\mathcal{T}(\mathcal{X})}$ such that $\tau_{\mathcal{X}}\circ w={\rm id}_{\mathcal{X}}$. In local coordinates, a vector field reads $w=w^i(x)\partial_{x^i}$. Likewise, a covector field or (differential) 1-form is a section of the cotangent bundle, i.e., a map $\omega\colon\mathcal{X}\rightarrow\mathcal{T}^\ast(\mathcal{X})$ which in local coordinates reads $\omega=\omega_i(x)\D x^i$. A codistribution on $\mathcal{X}$ of rank $k$ is a map which assigns to each $p\in\mathcal{X}$ a $k$-dimensional linear subspace $P_p\subset\mathcal{T}^\ast_p(\mathcal{X})$ of the cotangent space at~$p$. Then locally there exist $k$ covector fields $\omega^1,\dots,\omega^k$ such that $\omega^1_p,\dots,\omega^k_p$ form a basis for $P_p$. We say that the codistribution $P$ is (locally) spanned by the covector fields $\omega^1,\dots,\omega^k$, which form a (local) basis for $P$, i.e., $P=\Span{\omega^1,\dots,\omega^k}$, with the span over the ring $C^\infty(\mathcal{X})$ of smooth functions $h\colon\mathcal{X}\rightarrow\mathbb{R}$. The $k$-fold Lie derivative of a function $\varphi$ along a vector field $w$ is denoted by $\Lie_{w}^k\varphi$. By $\partial_xh$ we denote the $m\times n$ Jacobian matrix of $h=\bigl(h^1,\dots,h^m\bigr)$ with respect to $x=\bigl(x^1,\dots,x^n\bigr)$. The symbols $\subset$ and $\supset$ are used in the sense that they also include equality.
		We write $h_{[\alpha]}$ for the $\alpha$-th time derivative of a function $h$. When $h$ consists of several components, i.e., $h=\bigl(h^1,\dots,h^m\bigr)$, then \smash{$h_{[\alpha]}=\bigl(h^1_{[\alpha]},\dots,h^m_{[\alpha]}\bigr)$}. To keep expressions involving time derivatives of different orders short we use multi-indices. Let $A=\bigl(a^1,\dots,a^m\bigr)$ and $B=\bigl(b^1,\dots,b^m\bigr)$ be two multi-indices with $a^j\leq b^j$, $j=1,\dots,m$, which we abbreviate by~${A\leq B}$. Then
		\begin{gather*}
\begin{split}
& h_{[A]}=\bigl(h^1_{[a^1]},\dots,h^m_{[a^m]}\bigr),\qquad h_{[0,A]}=\bigl(h^1_{[0,a^1]},\dots,h^m_{[0,a^m]}\bigr),\\
& h_{[A,B]}=\bigl(h^1_{[a^1,b^1]},\dots,h^m_{[a^m,b^m]}\bigr),
\end{split}
		\end{gather*}
		where \smash{$h^j_{[a^j,b^j]}=\bigl(h^j_{[a^j]},\dots,h^j_{[b^j]}\bigr)$}. We define \smash{$h^j_{[\alpha^j,\beta^j]}$} to be empty when $\alpha^j>\beta^j$. Addition and subtraction of multi-indices is done componentwise and we define the addition and subtraction of a multi-index $A$ with an integer $c$ by $A\pm c=\bigl(a^1\pm c,\dots,a^m\pm c\bigr)$. Furthermore, we define~\smash{$|A|=\sum_{j=1}^{m}a^j$}. Multi-indices are sometimes also used in connection with the Lie derivative, for instance, if $f$ is a vector field, then $\Lie_f^Ah=\bigl(\Lie_f^{a^1}h^1,\dots,\Lie_f^{a^m}h^m\bigr)$. When $h$ consists of multiple blocks, the first subscript refers to the block and a second subscript in square brackets is used for denoting time derivatives. Consider, for instance, the function
			\begin{align*}
					h&=(h_1,h_2)=\bigl(\underbrace{h^1,\dots,h^{m_1}}_{h_{1}},\underbrace{h^{m_1+1},\dots,h^m}_{h_{2}}\bigr),
			\end{align*}
			whose components $h^j$ are grouped into the two blocks $h_1$ and $h_2$. Then a second subscript in square brackets is used for denoting time derivatives, e.g.,
			\begin{align*}
				h_{[\alpha]}
&=(h_{1,[\alpha]},h_{2,[\alpha]})=\bigl(\underbrace{h^1_{[\alpha]},\dots,h^{m_1}_{[\alpha]}}_{h_{1,{[\alpha]}}},\underbrace{h^{m_1+1}_{[\alpha]},\dots,h^m_{[\alpha]}}_{h_{2,{[\alpha]}}}\bigr)
		\end{align*}
		and \smash{$h_{1,[A_1]}=\bigl(h^1_{1,[a_1^1]},\dots,h^{m_1}_{1,[a_1^{m_1}]}\bigr)$} with some multi-index $A_1=\bigl(a_1^1,\dots,a_1^{m_1}\bigr)$.
		\begin{Example}
			Consider the tuple $h=\bigl(h^1,h^2\bigr)$, the integer $c=2$ and the multi-indices ${A=(1,3)}$, $B=(2,3)$. We then have
$h_{[c]}=\bigl(h^1_{[2]},h^2_{[2]}\bigr)$, $h_{[A]}=\bigl(h^1_{[1]},h^2_{[3]}\bigr)$, $h_{[0,A]}= \bigl(h^1,h^1_{[1]},h^2,\allowbreak h^2_{[1]}, h^2_{[2]},h^2_{[3]}\bigr)$, $h_{[A,B]}=\bigl(h^1_{[1]},h^1_{[2]},h^2_{[3]}\bigr)$,
 as well as $h_{[A+c]}=\bigl(h^1_{[3]},h^2_{[5]}\bigr)$ and $h_{[0,A-c]}=\bigl(h^2,h^2_{[1]}\bigr)$.
		\end{Example}
		
	\subsection[Geometric framework, flatness and quasi-static feedback of classical states]{Geometric framework, flatness\\ and quasi-static feedback of classical states}\label{sec:geometric_framework}
		Throughout this contribution, we use a finite-dimensional differential-geometric framework like, e.g., in \cite{KolarSchoberlSchlacher:2016-3}. All our definitions and results are generically local, that is, they hold locally in open and dense subsets. Thus, we can assume that all considered Jacobian matrices and (co-)distributions have locally constant rank. In order to compute time derivatives of functions of the system variables along trajectories of \eqref{eq:intro_nlsys}, a manifold $\XxUlu$ with coordinates~${(x,u,u_{[1]},\dots,u_{[l_u]})}$ is introduced, where $u_{[\alpha]}$ denotes the $\alpha$-th time derivative of the input $u$ and $l_u$ is some large enough but finite integer. The time derivative of a function $h(x,u,u_{[1]},\dots,u_{[l_u-1]})$, which does not explicitly depend on $u_{[l_u]}$, is then given by the Lie derivative $\Lie_{f_u}h$ along the vector field\footnote{To be precise, we have $\tfrac{\D}{\D t}(h\circ(x(t),u(t),u_{[1]}(t),\dots))=(\Lie_{f_u}h)\circ(x(t),u(t),u_{[1]}(t),\dots)$. Therefore, we identify the Lie derivative of a function on $\XxUlu$ along $f_u$ with its time derivative.}
		\begin{align}\label{eq:fu}
			\begin{aligned}
				f_u&=\sum_{i=1}^{n}f^i(x,u)\partial_{x^i}+\sum_{j=1}^{m}\sum_{\alpha=0}^{l_u-1}u^j_{[\alpha+1]}\partial_{u^j_{[\alpha]}}.
			\end{aligned}
		\end{align}
		In the remainder of the paper, we assume that $l_u$ is chosen large enough such that $f_u$ acts as time derivative on all functions considered. An appropriate bound on $l_u$ for the computations in this paper is given below Definition \ref{def:flatness}. Within this differential-geometric framework, flatness can be defined as follows.
		\begin{Definition}\label{def:flatness}
			System \eqref{eq:intro_nlsys} is called flat if there exists an $m$-tuple of smooth functions
			\begin{align}\label{eq:y}
					y^j&=\varphi^j(x,u,u_{[1]},\dots,u_{[q]}),\qquad j =1,\dots,m,
			\end{align}
			defined on $\XxUlu$ and smooth functions $F_x^i$ and $F_u^j$ such that locally
			\begin{align}
					x^i&=F_x^i\bigl(\varphi,\Lie_{f_u}\varphi,\dots,\Lie_{f_u}^{R-1}\varphi\bigr),\qquad i=1,\dots,n,\nonumber\\
					u^j&=F_u^j\bigl(\varphi,\Lie_{f_u}\varphi,\dots,\Lie_{f_u}^R\varphi\bigr),\qquad j=1,\dots,m,\label{eq:parameterization}
			\end{align}
			with some multi-index $R=\bigl(r^1,\dots,r^m\bigr)$. The $m$-tuple \eqref{eq:y} is called a flat output.
		\end{Definition}
		The highest time derivative of each component $y^j$ which may be needed throughout the derivations in this paper is given by $r^j$, where $r^j$ are the integers which form the multi-index $R=\bigl(r^1,\dots,r^m\bigr)$ in Definition \ref{def:flatness}. Therefore, a bound on $l_u$ is given by $r+q$ with ${r=\max\big\{r^1,\dots,r^m\big\}}$. A bound on $r$ in terms of the number of state variables $n$, the number of inputs $m$ and the highest order $q$ of time derivatives of the input occurring in \eqref{eq:y} can be found in \cite{KolarSchoberlSchlacher:2016-3}, namely $r\leq n+(m-1)q$. Hence, in total we have $n+(m-1)q+q=n+mq$ as a bound on $l_u$.
		
		We call a flat output $y=\varphi(x,u)$, which may depend on the input $u$ but not on time derivatives of the latter, an $(x,u)$-flat output and a system possessing such a flat output an $(x,u)$-flat system. It should be stressed that we do not require an $(x,u)$-flat output to explicitly depend on $u$. Therefore, every $x$-flat output $y=\varphi(x)$, resp.\ $x$-flat system, is also an $(x,u)$-flat output, resp.\ $(x,u)$-flat system (whereas the converse is obviously wrong). Our main results are derived for $(x,u)$-flat outputs and thus also apply to $x$-flat outputs.
		
		By taking the exterior derivative of expressions \eqref{eq:parameterization}, we find that flatness implies
		\begin{align*}
			\D x\in\operatorname{span}\big\{\D\varphi,\D\Lie_{f_u}\varphi,\dots,\D\Lie_{f_u}^{R-1}\varphi\big\}\qquad\text{and}\qquad\D 	u\in\operatorname{span}\big\{\D\varphi,\D\Lie_{f_u}\varphi,\dots,\D\Lie_{f_u}^R\varphi\big\},
		\end{align*}
		where throughout this section, span denotes the span over smooth functions on $\XxUlu$, i.e., smooth functions of $(x,u,u_{[1]},\dots)$. Under the assumption that ranks are constant (which we assume throughout), the converse is also true. This is a consequence of the following more general result relating the functional dependence of functions and the linear dependence of their differentials.
		\begin{Lemma}\label{lem:functional_independence}
			Consider a set of smooth functions $g^1,\dots,g^k$ as well as another smooth function~$h$ which are all defined on the same manifold. The following conditions are equivalent:
			\begin{enumerate}\itemsep=0pt
				\item[$(1)$] Locally $\D h\in\operatorname{span}\big\{\D g^1,\dots,\D g^k\big\}$.
				\item[$(2)$] There exists a smooth function $\psi\colon\mathbb{R}^k\rightarrow\mathbb{R}$ such that locally $h=\psi\bigl(g^1,\dots,g^k\bigr)$ holds identically.
			\end{enumerate}
			Furthermore, if the differentials $\D g^1,\dots,\D g^k$ are linearly independent, then the function $\psi$ is unique.
		\end{Lemma}
		The proof is straightforward and follows immediately by introducing a maximal number of independent functions of the set $g^1,\dots,g^k$ as local coordinates on the underlying manifold, see, e.g., \cite{Kolar:2017}. Another well-known important implication of Definition \ref{def:flatness} is that the differentials~\smash{$\D\varphi,\D\Lie_{f_u}\varphi,\dots,\D\Lie_{f_u}^\beta\varphi$} of derivatives of a flat output up to an arbitrary order $\beta$ are linearly independent.\footnote{The integer $l_u$ needs to be chosen large enough such that the Lie derivative $\Lie_{f_u}^\beta\varphi$ indeed yields the correct expression for the $\beta$-th time derivative of the functions $\varphi$. The linear independence of the differentials~\smash{$\D\varphi,\D\Lie_{f_u}\varphi,\dots,\D\Lie_{f_u}^\beta\varphi$} implies that the time evolution of the flat output is not constrained by any autonomous differential equation $\chi(y,y_{[1]},\dots,y_{[\beta]})=0$. We refer to the latter property as differential independence of the components of a flat output.} This in turn implies that there exists a unique minimal multi-index~${R=\bigl(r^1,\dots,r^m\bigr)}$, where $r^j$ is the order of the highest derivative of $y^j$ needed for expressing~$x$ and~$u$. Because of Lemma \ref{lem:functional_independence}, it furthermore follows that the corresponding maps $F_x$ and $F_u$ in \eqref{eq:parameterization} are unique. From now on, $R$ always denotes this unique minimal multi-index.
		
		Finally, let us give definitions for the notions of endogenous feedback and quasi-static feedback of a classical state.
		\begin{Definition}\label{def:endogenous}
				Consider system \eqref{eq:intro_nlsys} and a feedback of the form
				\begin{align}\label{eq:endogenous feedback}
						\dot{z}&=g(x,z,v,v_{[1]},\dots),\qquad
						u=\alpha(x,z,v,v_{[1]},\dots),
				\end{align}
				which contains the dynamic feedback \eqref{eq:intro_dynamic_feedback} and the quasi-static feedback \eqref{eq:intro_quasi_static_feedback_x} as special cases. The feedback \eqref{eq:endogenous feedback} is called endogenous if its state $z$ and the new input $v=\bigl(v^1,\dots,v^m\bigr)$ can be expressed as functions of $x$, $u$ and time derivatives of $u$, see also \cite[Definition 2.3]{Martin:1992}.
		\end{Definition}
		\begin{Definition}\label{def:quasistatic x}
			Consider system \eqref{eq:intro_nlsys} with state $x$, control input $u$ and an invertible transformation of the form
			\begin{alignat}{3}
				&	u=\alpha(x,v,v_{[1]},\dots),	\qquad && v=\hat{\alpha}(x,u,u_{[1]},\dots),&\nonumber\\
					&u_{[1]}=\alpha_1(x,v,v_{[1]},\dots),\qquad&& v_{[1]}=\hat{\alpha}_1(x,u,u_{[1]},\dots),&\nonumber\\
				&	u_{[2]}=\alpha_2(x,v,v_{[1]},\dots), \qquad&& v_{[2]}=\hat{\alpha}_2(x,u,u_{[1]},\dots),&\nonumber\\
					&\vdotswithin{=} \qquad &&  \vdotswithin{=}&
\label{eq:quasi_static_feedback_x}
			\end{alignat}
			such that the differentials $\D x,\D v,\D v_{[1]},\dots$ are linearly independent. We refer to $u=\alpha(x,v,\allowbreak v_{[1]},\dots)$ in \eqref{eq:quasi_static_feedback_x} as a quasi-static feedback of the state $x$.\footnote{Note that \eqref{eq:quasi_static_feedback_x} implies $\Dim{v}=\Dim{u}$.}
		\end{Definition}

		\begin{Remark}
			A quasi-static feedback of a classical state as in Definition \ref{def:quasistatic x} is a special case of a quasi-static feedback as considered, e.g., in \cite{DelaleauRudolph:1995,DelaleauRudolph:1998}. The major difference is that our point of departure is always a system of the form \eqref{eq:intro_nlsys} in classical state representation, whereas in \cite{DelaleauRudolph:1995,DelaleauRudolph:1998} generalized state representations $\dot{\tilde{x}}=f(\tilde{x},u,u_{[1]},\dots)$ are considered. The transformation \eqref{eq:quasi_static_feedback_x} would be the same in the case of generalized state representations with $x$ replaced by $\tilde{x}$.
		\end{Remark}
		
		\section{Exact feedback linearization of flat systems}\label{sec:exact_lin}
		In this section, we discuss the exact feedback linearization of a~system \eqref{eq:intro_nlsys} with respect to a~given flat output \eqref{eq:y} by a~suitable feedback. The feedback-modified system shall possess a~linear input-output behavior of the form \smash{$y^j_{[a^j]}=v^j$, $j=1,\dots,m$} -- with suitable integers $a^j$ -- between a~newly introduced input $v$ and the flat output $y$. Such an exact feedback linearization is used, e.g., as an intermediate step in the design of a flatness-based tracking control. For the design of a tracking control, it is desirable to choose the orders $a^j$ of the time derivatives of the components of the flat output which are used as new (closed-loop) input $v$ as low as possible to obtain tracking error dynamics of minimal order. In Theorem \ref{thm:linerizing_feedback} below, we provide easily verifiable conditions which assure that a selection of time derivatives of a flat output can be introduced as new (closed-loop) input $v$. The construction of a corresponding feedback which introduces these time derivatives as new input is similar to the construction of the classical linearizing endogenous dynamic feedback \eqref{eq:intro_dynamic_feedback} as proposed, e.g., in \cite{Chetverikov:2007,FliessLevineMartinRouchon:1999,Kolar:2017}. So let us first briefly recall its construction, which is as follows.
		
		We have to extend the flat parameterization $F_x$ of the state of the system to a diffeomorphism, i.e., we have to choose $|R|-n$ functions $z^l=F_z^l(\varphi_{[0,R-1]})$ such that the map $(F_x,F_z)\colon\mathbb{R}^{|R|}\rightarrow\mathbb{R}^{|R|}$ is a diffeomorphism. We then have
		\begin{align}
				&x^i=F_x^i(\varphi_{[0,R-1]}),\qquad i=1,\dots,n,\nonumber\\
				& z^l=F_z^l(\varphi_{[0,R-1]}),\qquad l=1,\dots,|R|-n,\label{eq:Fxz}
		\end{align}
		which allows us to express the functions \smash{$\varphi_{[0,R-1]}$} as functions of $x$ and $z$, and by setting~${v=\varphi_{[R]}}$, all the functions $\varphi_{[0,R]}$ can be expressed as functions of $x$, $z$ and $v$, i.e., there is a diffeo\-mor\-phism \smash{$\Psi\colon \mathbb{R}^{|R|+m}\rightarrow\mathbb{R}^{|R|+m}$} such that $\varphi_{[0,R]}=\Psi(x,z,v)$. The time derivatives $\smash{\dot{z}^l}=\allowbreak\smash{\Lie_{f_u}\bigl(F_z^l(\varphi_{[0,R-1]})\bigr)}$ are obviously also functions of $\varphi_{[0,R]}=\Psi(x,z,v)$, i.e., with
\[g^l(x,z,v):=\Lie_{f_u}\bigl(F_z^l(\varphi_{[0,R-1]})\bigr)\circ\Psi(x,z,v),\qquad l=1,\dots,|R|-n,
\]
 we have $\dot{z}^l=g^l(x,z,v)$. Now consider the dynamic feedback
		\begin{align}\label{eq:endogenous_dynamic_feedback_2}
				\dot{z}&=g(x,z,v),\qquad
				u=F_u\circ\Psi(x,z,v),
		\end{align}
		constructed from the functions $g$, the flat parameterization $F_u$ of the input and the map $\Psi$. Applying this feedback to system \eqref{eq:intro_nlsys} yields the closed-loop system
		\begin{align}\label{eq:endogenous_dynamic_feedback_closed_loop_2}
				\dot{x}&=f(x,F_u\circ\Psi(x,z,v)),\qquad
				\dot{z}=g(x,z,v)
		\end{align}
		with the state $(x,z)$ and the new input $v$. It can then be show that the transformation $y_{[0,R]}=\Psi(x,z,v)$ puts the closed-loop system into the Brunovsk\'y normal form
		\begin{align}\label{eq:closed_loop_system_Brunovsky_dynamic}
				\dot{y}_{[0,R-1]}&=y_{[1,R]}.
		\end{align}
		The top variables $y^j=\Psi^j(x,z,v)$, $j=1,\dots,m$ of the integrator chains in \eqref{eq:closed_loop_system_Brunovsky_dynamic} obviously form a~linearizing output of the closed-loop system. These top variables are in fact just the components of the flat output $\varphi(x,u,u_{[1]},\dots)$ of the original system expressed in terms of the state $(x,z)$ and the input $v$ of the closed-loop system \eqref{eq:endogenous_dynamic_feedback_closed_loop_2}. The closed-loop system thus has the linear input-output behavior $y_{[R]}=v$ between the new input $v$ and the flat output~${y^j=\varphi^j(x,u,u_{[1]},\dots)=\Psi^j(x,z,v)}$. In other words, feedback \eqref{eq:endogenous_dynamic_feedback_2} introduces the highest time derivatives $\varphi_{[R]}$ which are present in \eqref{eq:parameterization} as the new input $v$.
		
		As already mentioned, for the design of a tracking control, it is desirable to choose the orders $a^j$ of the time derivatives of the components of the flat output which are used as new (closed-loop) input $v$ as low as possible to obtain tracking error dynamics of minimal order. The following theorem provides easily verifiable conditions which assure that a given selection of time derivatives of a flat output can be introduced as new (closed-loop) input $v$. The multi-index $R$ in this theorem again refers to the minimal multi-index such that \eqref{eq:parameterization} in Definition \ref{def:flatness} holds.
		\begin{Theorem}\label{thm:linerizing_feedback}
			Consider system \eqref{eq:intro_nlsys} with flat output $y=\varphi(x,u,u_{[1]},\dots,u_{[q]})$. For any $m$-tuple $y_{[A]}$ with $A\leq R$, satisfying the property that the differentials $\D x,\D\varphi_{[A]},\D\varphi_{[A+1]},\dots,\D\varphi_{[R-1]}$ are linearly independent, there exists an endogenous feedback of the form
				\begin{align}\label{eq:endogenous_dynamic_quasistatic_feedback}
						\dot{z}&=g(x,z,v,v_{[1]},\dots),\qquad
						u=\alpha(x,z,v,v_{[1]},\dots)
				\end{align}
				with $\Dim{z}=|A|-n$, such that the closed-loop system has the input-output behavior $y_{[A]}=v$. In the case $|A|=n$, $z$ is empty and the feedback \eqref{eq:endogenous_dynamic_quasistatic_feedback} reduces to $u=\alpha(x,v,v_{[1]},\dots)$, i.e., a~quasi-static feedback of the state $x$.
		\end{Theorem}

		The proof of this theorem is done constructively and thereby provides a systematic method for constructing a feedback which actually introduces the time derivatives $\varphi_{[A]}$ of a~flat output as new input. The condition on the linear independence of the differentials $\D x,\D\varphi_{[A]},\D\varphi_{[A+1]},\allowbreak\dots,\D\varphi_{[R-1]}$\footnote{Note that due to the linear independence of the differentials of time derivatives of a flat output up to an arbitrary order and $\D x\in\Span{\D\varphi_{[0,R-1]}}$, the linear independence of $\D x,\D\varphi_{[A]},\D\varphi_{[A+1]},\dots,\D\varphi_{[R-1]}$ actually implies linear independence of the differentials $\D x,\D\varphi_{[A]},\D\varphi_{[A+1]},\dots,\D\varphi_{[\beta]}$ up to an arbitrary order $\beta$.} roughly speaking assures that trajectories of $\varphi_{[A]}$ (i.e., the derivatives of the flat output intended as new inputs) are not restricted by the current state $x$ of the system. This is analogous to the independence of the original control input $u$ and the state $x$ (the differentials $\D x,\D u,\D u_{[1]},\dots$ are obviously linearly independent). Note that the time derivatives $\varphi_{[R]}$ introduced as new input by means of the endogenous dynamic feedback discussed above always satisfy this independence condition. Furthermore, it can be shown that the condition can only be satisfied for multi-indices $A$ with $|A|\geq n$. An obvious example which does not meet this condition are the components of an $x$-flat output.
		\begin{proof}[Proof of Theorem \ref{thm:linerizing_feedback}]
			The construction of the feedback is similar to the construction of the classical endogenous dynamic feedback which introduces the time derivatives $\varphi_{[R]}$ of a~flat output as new input and which we recalled above. However, instead of extending the flat parameterization $x=F_x(\varphi_{[0,R-1]})$ of the state to a diffeomorphism with $|R|-n$ arbitrary functions~\smash{$z^l\!=\!F_z^l(\varphi_{[0,R-1]})$} (see \eqref{eq:Fxz}), here, a part of the functions used for extending~${x\!=\! F_x(\varphi_{[0,R-1]})}$ to a diffeomorphism will be given by the consecutive time derivatives $\varphi_{[A,R-1]}$.
			
			By assumption, the differentials $\D x$, $\D\varphi_{[A,R-1]}$ are linearly independent. Because of that and~${x=F_x(\varphi_{[0,R-1]})}$, it follows that the map $\varphi_{[0,R-1]}\mapsto(F_x(\varphi_{[0,R-1]}),\varphi_{[A,R-1]})$ is a submersion. Therefore, we can chose $|A|-n$ further functions \smash{$z^l=F_z^l(\varphi_{[0,R-1]})$}, $l=1,\dots,|A|-n$ such that the map $\varphi_{[0,R-1]}\mapsto(F_x(\varphi_{[0,R-1]}),F_z(\varphi_{[0,R-1]}),\varphi_{[A,R-1]})$ is a diffeomorphism. The map $\varphi_{[0,R]}\mapsto(F_x(\varphi_{[0,R-1]}),F_z(\varphi_{[0,R-1]}),\varphi_{[A,R]})$ is then obviously also a diffeomorphism. Via the inverse of this diffeomorphism, all the functions $\varphi_{[0,R]}$ can be expressed as functions of $x$, $z=F_z(\varphi_{[0,R-1]})$ and $v_{[0,R-A]}:=\varphi_{[A,R]}$, i.e., there exists a diffeomorphism $\Psi$ such that
			\begin{align*}
					\varphi_{[0,R]}&=\Psi(x,z,v_{[0,R-A]}).
			\end{align*}
			The time derivatives \smash{$\dot{z}^l=\Lie_{f_u}\bigl(F_z^l(\varphi_{[0,R-1]})\bigr)$} are obviously also functions of $\varphi_{[0,R]}=\Psi(x,z,\allowbreak v_{[0,R-A]})$, i.e., with \smash{$g^l(x,z,v_{[0,R-A]}):=\Lie_{f_u}\bigl(F_z^l(\varphi_{[0,R-1]})\bigr)\circ\Psi(x,z,v_{[0,R-A]})$}, $l=1,\dots,|A|-n$, we have \smash{$\dot{z}^l=g^l(x,z,v_{[0,R-A]})$}.
			Now consider the dynamic feedback
			\begin{align}\label{eq:dynamic_feedback_v2}
					\dot{z}&=g(x,z,v_{[0,R-A]}),\qquad
					\dot{v}_{[0,R-A-1]}=v_{[1,R-A]},\qquad
					u=F_u\circ\Psi(x,z,v_{[0,R-A]})
			\end{align}
			constructed from the functions $g$, the flat parameterization $F_u$ of the input and the map $\Psi$. Applying this feedback to system \eqref{eq:intro_nlsys} yields the closed-loop system
			\begin{align}\label{eq:closed_loop_system_dynamic_feedback_v2}
					\dot{x}&=f(x,F_u\circ\Psi(x,z,v_{[0,R-A]})),\qquad
					\dot{z}=g(x,z,v_{[0,R-A]}),\qquad
					\dot{v}_{[0,R-A-1]}=v_{[1,R-A]}
			\end{align}
			with the state $(x,z,v_{[0,R-A-1]})$ and the new input $v_{[R-A]}$. Recall that feedback \eqref{eq:dynamic_feedback_v2} is just a~special case of the endogenous dynamic feedback \eqref{eq:endogenous_dynamic_feedback_2} considered above, the difference being that a part of the functions used for extending the map $F_x$ to a diffeomorphism is given by the consecutive time derivatives $\varphi_{[A,R-1]}$ instead of arbitrary functions. It is thus immediate that the closed-loop system \eqref{eq:closed_loop_system_dynamic_feedback_v2} again has a linear input-output behavior between the new input~$v_{[R-A]}$ and the flat output $y$.
			
			Now consider the feedback
			\begin{align*}
					\dot{z}&=g(x,z,v_{[0,R-A]}),\qquad
					u=F_u\circ\Psi(x,z,v_{[0,R-A]}),
			\end{align*}
			obtained from \eqref{eq:dynamic_feedback_v2} by omitting the integrator chains $\dot{v}_{[0,R-A-1]}=v_{[1,R-A]}$. This feedback yields the closed-loop system
			\begin{align}\label{eq:closed_loop_system_dynamic_quasistatic_feedback_v2}
					\dot{x}&=f(x,F_u\circ\Psi(x,z,v_{[0,R-A]})),\qquad
					\dot{z}=g(x,z,v_{[0,R-A]})
			\end{align}
			with the new input $v$ (instead of $v_{[R-A]}$). Computing the time derivative of a function along trajectories of \eqref{eq:closed_loop_system_dynamic_quasistatic_feedback_v2} or along trajectories of \eqref{eq:closed_loop_system_dynamic_feedback_v2} yields exactly the same result. For the closed-loop system \eqref{eq:closed_loop_system_dynamic_quasistatic_feedback_v2}, we thus have the linear input-output behavior $y_{[A]}=v$.
		\end{proof}

		\begin{Remark}
			Theorem \ref{thm:linerizing_feedback} analogously applies to multi-indices $A$ which do not meet ${A\leq R}$, only minor modifications of the proof would be required in this case. However, since the choice~${A=R}$ is always possible, and in regard of a subsequently designed tracking control a~choice for~$A$ which is larger than necessary does not seem to be desirable, we for simplicity assume $A\leq R$.
		\end{Remark}

		For the design of a flatness-based tracking control, the special case $|A|=n$ is particularly interesting, since it ensures a tracking error dynamics of minimal order. Furthermore, it has the advantage that according to Theorem \ref{thm:linerizing_feedback} the required feedback is a quasi-static one and no controller states need to be initialized. Even though it has been shown in \cite{DelaleauRudolph:1998} that every flat system with a flat output of the general form \eqref{eq:y} can be exactly linearized by a quasi-static feedback of a generalized state, the question whether this is also always possible by a~quasi-static feedback of a classical state is still open. For the practically most relevant subclass of~$(x,u)$-flat systems with flat outputs of the form \eqref{eq:xu_flat_output}, however, the existence of a linearizing quasi-static feedback of a classical state can be deduced by applying results derived in \cite{DelaleauPereiradaSilva:1998-1,DelaleauPereiradaSilva:1998-2} in a differential-algebraic setting. In the following subsection, we provide a self-contained proof as well as a~systematic construction of such a quasi-static feedback in our geometric framework.
		\subsection[Exact feedback linearization of (x,u)-flat systems]{Exact feedback linearization of $\boldsymbol{(x,u)}$-flat systems}\label{subsec:exact_lin_xu}
		
		In this section, it is shown in our finite-dimensional geometric framework that every system \eqref{eq:intro_nlsys} with an $(x,u)$-flat output
		\begin{align}\label{eq:xu_flat_output_procedure}
			\begin{aligned}
				y&=\varphi(x,u)
			\end{aligned}
		\end{align}
		can be exactly feedback linearized with respect to this output by a quasi-static feedback of its original state $x$. Because of Theorem \ref{thm:linerizing_feedback}, we only have to show that there exists a multi-index $A\leq R$ with $|A|=n$ such that the differentials $\D x,\D\varphi_{[A]},\D\varphi_{[A+1]},\dots,\D\varphi_{[R-1]}$ are linearly independent. From now on it is referred to such a special multi-index by $\kappa=\bigl(\kappa^1,\dots,\kappa^m\bigr)$ in order to distinguish it from the general case.
		
		In the following, a procedure for systematically constructing such a multi-index is proposed. In each step of the procedure, certain time derivatives of the flat output are introduced as new coordinates on $\XxUlu$, such that finally the coordinates $(x,u,u_{[1]},\dots)$ are replaced by~${(x,v,v_{[1]},\dots)}$ with $v=y_{[\kappa]}$. The multi-index $\kappa$ is of course not known a priori, it is constructed successively and thus also the new coordinates $v$ are introduced successively. Roughly speaking, we differentiate each component of the flat output until it depends explicitly on the input $u$. By means of a coordinate transformation, we then replace as many components of the input $u$ as possible by these derivatives. In the next step, each of the remaining components of the flat output is further differentiated until again an explicit dependence on the remaining components of the original input $u$ occurs, and again as many of its components as possible are replaced by these derivatives. This procedure is continued until all components of the original input $u$ have been replaced by time derivatives of the flat output. For the proposed procedure, it is crucial that the ranks of the occurring Jacobian matrices are locally constant, which we assume throughout.
		
		\textit{Step $1$.} Define the multi-index $K_1=\bigl(k_1^1,\dots,k_1^m\bigr)$ such that
		\begin{align*}
				\Lie_{f_u}^{k_1^j-1}\varphi^j&=\varphi^j_{[k_1^j-1]}(x),\qquad
				\Lie_{f_u}^{k_1^j}\varphi^j=\varphi^j_{[k_1^j]}(x,u),
		\end{align*}
		i.e., \smash{$k_1^j$} denotes the relative degree of the component $\varphi^j$ of the flat output. Note that for an~$(x,u)$-flat output where $y^j=\varphi^j(x,u)$ explicitly depends on $u$, we have $k_1^j=0$. Introduce~${m_1=\Rank{\partial_u\varphi_{[K_1]}}}$ of the $m$ functions $\varphi_{[K_1]}$ and their time derivatives as new coordinates. By reordering the components of the flat output we can always achieve that $\Rank{\partial_u\varphi_{1,[\kappa_1]}}=m_1$, where $\varphi_1=\bigl(\varphi^1,\dots,\varphi^{m_1}\bigr)$ and $\kappa_1=\bigl(k_1^1,\dots,k_1^{m_1}\bigr)$ consists of the first $m_1$ integers in $K_1$, enabling us to apply the coordinate transformation
		\begin{alignat}{3}
				&v_1=\varphi_{1,[\kappa_1]}(x,u),\qquad&&
				u_{{\rm rest}_1}=\bigl(u^{m_1+1},\dots,u^m\bigr),&\nonumber\\
				&v_{1,[1]}=\varphi_{1,[\kappa_1+1]}(x,u,u_{[1]}),\qquad&&
				u_{{\rm rest}_1,[1]}=\bigl(u^{m_1+1}_{[1]},\dots,u^m_{[1]}\bigr),&\nonumber\\
				&v_{1,[2]}=\varphi_{1,[\kappa_1+2]}(x,u,u_{[1]},u_{[2]}),\qquad&&
				u_{{\rm rest}_1,[2]}=\bigl(u^{m_1+1}_{[2]},\dots,u^m_{[2]}\bigr),&\nonumber\\
				&\vdotswithin{=}\qquad && \vdotswithin{=} & \label{eq:trnas_1}
		\end{alignat}
		That is, we replace the inputs $u^1,\dots,u^{m_1}$ and their time derivatives by $v_1$ and its time derivatives (which may require a renumbering of the inputs).\footnote{Introducing not only $v_1$ but also its time derivatives as new coordinates is crucial for preserving the simple structure of the vector field \eqref{eq:fu}. Concerning the regularity of transformation \eqref{eq:trnas_1} it should be noted that $\Rank{\partial_u\varphi_{1,[\kappa_1]}}=m_1$ implies \smash{$\Rank{\partial_{u_{[\alpha]}}\varphi_{1,[\kappa_1+\alpha]}}=m_1$} for $\alpha\geq 1$. In the new coordinates, the vector field \eqref{eq:fu} has the form
\[
f_u=f^i\circ\hat{\Phi} \partial_{x^i}+\sum_{\alpha=0}^{l_u-1}\bigl(v_{1,[\alpha+1]}^{j_1}\partial_{v_{1,[\alpha]}^{j_1}}
+u_{{\rm rest}_1,[\alpha+1]}^{j_2}\partial_{u_{{\rm rest}_1,[\alpha]}^{j_2}}\bigr)+\dots+\partial_{v_{1,[l_u]}}
\]
 with the inverse \smash{$\hat{\Phi}$} of the transformation \eqref{eq:trnas_1}. The in general non-zero components in the $\partial_{v_{1,[l_u]}}$-directions do not bother us as long as $l_u$ is chosen large enough.} In these coordinates, we have
		\begin{alignat*}{3}
				&y_{1,[0,\kappa_1-1]}=\varphi_{1,[0,\kappa_1-1]}(x),\qquad &&	y_{{\rm rest}_1,[0,K_{{\rm rest}_1}-1]}=\varphi_{{\rm rest}_1,[0,K_{{\rm rest}_1}-1]}(x),& \\
				&y_{1,[\kappa_1]}=v_1, \qquad && y_{{\rm rest}_1,[K_{{\rm rest}_1}]}=\varphi_{{\rm rest}_1,[K_{{\rm rest}_1}]}(x,v_1),&
		\end{alignat*}
		where
\[
y_1=\bigl(y^1,\dots,y^{m_1}\bigr),\qquad y_{{\rm rest}_1}=\bigl(y^{m_1+1},\dots,y^m\bigr), \qquad\varphi_{{\rm rest}_1,[\beta]}=\bigl(\varphi^{m_1+1}_{[\beta]},\dots,\varphi^m_{[\beta]}\bigr)\circ\hat{\Phi}
\]
 with the inverse \smash{$\hat{\Phi}$} of transformation \eqref{eq:trnas_1}, and $K_{{\rm rest}_1}\!\!=\!\!\bigl(k_1^{m_1+1},\dots,k_1^m\bigr)$. Note that $\varphi_{{\rm rest}_1,[K_{{\rm rest}_1}]}$ is independent of $u_{{\rm rest}_1}$, since otherwise $\Rank{\partial_u\varphi_{[K_1]}}$ would have been larger than $m_1$.
		
		\textit{Step $2$.} Define the multi-index $K_2=\bigl(k_2^1,\dots,k_2^{m-m_1}\bigr)$ such that
		\begin{align*}
				&\Lie_{f_u}^{k_2^j-1}\varphi_{{\rm rest}_1}^j=\varphi_{{\rm rest}_1,[k_2^j-1]}^j(x,v_1,v_{1,[1]},\dots),
\\
				&\Lie_{f_u}^{k_2^j}\varphi_{{\rm rest}_1}^j=\varphi_{{\rm rest}_1,[k_2^j]}^j(x,v_1,v_{1,[1]},\dots,u_{{\rm rest}_1}).
		\end{align*}
		Similar as before, we introduce $m_2=\Rank{\partial_{u_{{\rm rest}_1}}\varphi_{{\rm rest}_1,[K_2]}}$ of the $m-m_1$ functions $\varphi_{{\rm rest}_1,[K_2]}$ and their time derivatives as new coordinates. By reordering the components of the flat output belonging to $y_{{\rm rest}_1}$, we can always achieve that $\Rank{\partial_{u_{{\rm rest}_1}}\varphi_{2,[\kappa_2]}}=m_2$, where $\varphi_2=\bigl(\varphi_{{\rm rest}_1}^1,\dots,\varphi_{{\rm rest}_1}^{m_2}\bigr)$ and $\kappa_2=\bigl(k_2^1,\dots,k_2^{m_2}\bigr)$ consists of the first $m_2$ integers in $K_2$, enabling us to apply the coordinate transformation\footnote{The functions $\varphi_{2,[\kappa_2]}$ depend on time derivatives of $v_1$, which are only available up to the order $l_u$. Thus, some higher-order time derivatives of $u_{{\rm rest}_1}$ must be kept as coordinates on the finite-dimensional manifold $\XxUlu$ and cannot be replaced by time derivatives of $v_2$. However, this is no problem as long as $l_u$ is chosen large enough.}
		\begin{align}
&v_{2}=\varphi_{2,[\kappa_2]}(x,v_{1},v_{1,[1]},\dots,u_{{\rm rest}_1}),\nonumber\\
&u_{{\rm rest}_2}=\bigl(u^{m_1+m_2+1},\dots,u^m\bigr),\nonumber\\
&v_{2,[1]}=\varphi_{2,[\kappa_2+1]}(x,v_{1},v_{1,[1]},\dots,u_{{\rm rest}_1},u_{{\rm rest}_1,[1]}),\nonumber\\
&u_{{\rm rest}_2,[1]}=\bigl(u^{m_1+m_2+1}_{[1]},\dots,u^m_{[1]}\bigr),\nonumber\\
&v_{2,[2]}=\varphi_{2,{[\kappa_2+2]}}(x,v_{1},v_{1,[1]},\dots,u_{{\rm rest}_1},u_{{\rm rest}_1,[1]},u_{{\rm rest}_1,[2]}),\nonumber\\
&u_{{\rm rest}_2,[2]}=\bigl(u^{m_1+m_2+1}_{[2]},\dots,u^m_{[2]}\bigr),\nonumber\\
&\vdotswithin{=}\label{eq:trnas_2}
\end{align}
That is, we replace the inputs $u^{m_1+1},\dots,u^{m_1+m_2}$ and their time derivatives by $v_2$ and its time derivatives (which may again require a renumbering of the inputs belonging to $u_{{\rm rest}_1}$). In these coordinates we have
		\begin{gather*}
				y_{1,[0,\kappa_1-1]}=\varphi_{1,[0,\kappa_1-1]}(x),\\
				y_{1,[\kappa_1]}=v_{1},\\
				y_{2,[0,\kappa_2-1]}=\varphi_{2,[0,\kappa_2-1]}(x,v_{1},v_{1,[1]},\dots),\\
				y_{2,[\kappa_2]}=v_{2},\\
				y_{{\rm rest}_2,[0,K_{{\rm rest}_2}-1]}=\varphi_{{\rm rest}_2,[0,K_{{\rm rest}_2}-1]}(x,v_{1},v_{1,[1]},\dots),\\
				y_{{\rm rest}_2,[K_{{\rm rest}_2}]}=\varphi_{{\rm rest}_2,[K_{{\rm rest}_2}]}(x,v_{1},v_{1,[1]},\dots,v_{2}),
		\end{gather*}
		where
		\begin{gather*}
y_2=\bigl(y^{m_1\!+1},\dots,y^{m_1\!+m_2}\bigr),\qquad y_{{\rm rest}_2}=\bigl(y^{m_1+m_2+1},\dots,y^m\bigr), \\ \varphi_{{\rm rest}_2,[\beta]=\bigl(\varphi_{{\rm rest}_1,[\beta]}^{m_2+1},\dots,\varphi_{{\rm rest}_1,[\beta]}^{m-m_1}\bigr)}\allowbreak\circ\hat{\Phi}
\end{gather*}
 with the inverse \smash{$\hat{\Phi}$} of transformation \eqref{eq:trnas_2} and $K_{{\rm rest}_2}=\bigl(k_2^{m_2+1},\dots,k_2^{m-m_1}\bigr)$.
		
		\textit{Step $i$.}	In the $i$-th step, we are concerned with the dependence of the time derivatives of the functions $\varphi_{{\rm rest}_{i-1}}\!=\bigl(\varphi_{{\rm rest}_{i-2}}^{m_{i-1}+1},\dots,\varphi_{{\rm rest}_{i-2}}^{m-m_1-\dots-m_{i-2}}\bigr)$ on the inputs $u_{{\rm rest}_{i-1}}\!= \bigl(u^{m_1+\dots+m_{i-1}+1},\dots,\allowbreak u^m\bigr)$. Define the multi-index \smash{$K_i=\bigl(k_i^1,\dots,k_i^{m-m_1-\dots-m_{i-1}}\bigr)$} such that
		\begin{align*}
				&\Lie_{f_u}^{k_i^j-1}\varphi_{{\rm rest}_{i-1}}^j=\varphi_{{\rm rest}_{i-1},[k_i^j-1]}^j(x,v_{1},v_{1,[1]},\dots,v_{i-1},v_{i-1,[1]},\dots),\\
				&\Lie_{f_u}^{k_i^j}\varphi_{{\rm rest}_{i-1}}^j=\varphi_{{\rm rest}_{i-1},[k_i^j]}^j(x,v_{1},v_{1,[1]},\dots,v_{i-1},v_{i-1,[1]},\dots,u_{{\rm rest}_{i-1}}).
		\end{align*}
		Introduce $m_i=\Rank{\partial_{u_{{\rm rest}_{i-1}}}\varphi_{{\rm rest}_{i-1},[K_i]}}$ of the $m-m_1-\dots-m_{i-1}$ functions $\varphi_{{\rm rest}_{i-1},[K_i]}$ and their time derivatives as new coordinates. By reordering the components of the flat output belonging to $\varphi_{{\rm rest}_{i-1}}$, we can always achieve that \smash{$\Rank{\partial_{u_{{\rm rest}_{i-1}}}\varphi_{i,[\kappa_i]}}\!=\!m_i$}, where $\smash{\varphi_i\!=\!\bigl(\varphi_{{\rm rest}_{i-1}}^1,\dots,}\allowbreak\smash{\varphi_{{\rm rest}_{i-1}}^{m_i}\bigr)}$ and $\kappa_i=\bigl(k_i^1,\dots,k_i^{m_i}\bigr)$ consists of the first $m_i$ integers in $K_i$, enabling us to apply the coordinate transformation
		\begin{gather}
				v_{i}=\varphi_{i,[\kappa_i]}(x,v_{1},v_{1,[1]},\dots,v_{i-1},v_{i-1,[1]},\dots,u_{{\rm rest}_{i-1}}),\nonumber\\
				u_{{\rm rest}_i}=\bigl(u^{m_1+\dots+m_i+1},\dots,u^m\bigr),\nonumber\\
				v_{i,[1]}=\varphi_{i,[\kappa_i+1]}(x,v_{1},v_{1,[1]},\dots,v_{i-1},v_{i-1,[1]},\dots,u_{{\rm rest}_{i-1}},u_{{\rm rest}_{i-1},[1]}),\nonumber\\
				u_{{\rm rest}_i,[1]}=\bigl(u^{m_1+\dots+m_i+1}_{[1]},\dots,u^m_{[1]}\bigr),\nonumber\\
				v_{i,[2]}=\varphi_{i,[\kappa_i+2]}(x,v_{1},v_{1,[1]},\dots,v_{i-1},v_{i-1,[1]},\dots,u_{{\rm rest}_{i-1}},u_{{\rm rest}_{i-1},[1]},u_{{\rm rest}_{i-1},[2]}),\nonumber\\
				u_{{\rm rest}_i,[2]}=\bigl(u^{m_1+\dots+m_i+1}_{[2]},\dots,u^m_{[2]}\bigr),\nonumber\\
				\vdotswithin{=}\label{eq:trnas_i}
		\end{gather}
		That is, we replace the inputs \smash{$u^{m_1+\dots+m_{i-1}+1},\dots,u^{m_1+\dots+m_i}$} and their time derivatives by $v_i$ and its time derivatives (which may require a renumbering of the inputs belonging to $u_{{\rm rest}_{i-1}}$). In these coordinates, we have
		\begin{gather*}
				y_{1,[0,\kappa_1-1]} =\varphi_{1,[0,\kappa_1-1]}(x),\\
				y_{1,[\kappa_1]} =v_{1},\\
				y_{2,[0,\kappa_2-1]} =\varphi_{1,[0,\kappa_2-1]}(x,v_{1},v_{1,[1]},\dots),\\
				y_{2,[\kappa_2]} =v_{2},\\
				 \vdotswithin{=}\\		
	y_{i,[0,\kappa_i-1]} =\varphi_{i,[0,\kappa_i-1]}(x,v_{1},v_{1,[1]},\dots,v_{i-1},v_{i-1,[1]},\dots),\\
				y_{i,[\kappa_i]} =v_{i},\\
			y_{{\rm rest}_i,[0,K_{{\rm rest}_i}-1]} =\varphi_{{\rm rest}_i,[0,K_{{\rm rest}_i}-1]}(x,v_{1},v_{1,[1]},\dots,v_{i-1},v_{i-1,[1]},\dots),\\
				y_{{\rm rest}_i,[K_{{\rm rest}_i}]} =\varphi_{{\rm rest}_i,[K_{{\rm rest}_i}]}(x,v_{1},v_{1,[1]},\dots,v_{i-1},v_{i-1,[1]},\dots,v_{i}),
		\end{gather*}
		where
		\begin{gather*}
y_i=\bigl(y^{m_1+\dots+m_{i-1}+1},\dots,y^{m_1+\dots+m_i}\bigr),\qquad y_{{\rm rest}_i}=\bigl(y^{m_1+\dots+m_i+1},\dots,y^m\bigr), \\ \varphi_{{\rm rest}_i,[\beta]}=\bigl(\varphi_{{\rm rest}_{i-1},[\beta]}^{m_i+1},\dots,\varphi_{{\rm rest}_{i-1},[\beta]}^{m-m_1-\dots-m_{i-1}}\bigr)\circ\hat{\Phi}
\end{gather*}
 with the inverse \smash{$\hat{\Phi}$} of transformation~\eqref{eq:trnas_i} and $\smash{K_{{\rm rest}_i}=\bigl(k_i^{m_i+1},}\allowbreak\dots,\smash{k_i^{m-m_1-\dots-m_{i-1}}\bigr)}$.

\textit{Last step.} The procedure terminates when in some step, let us call it the $s$-th step, the Jacobian matrix $\partial_{u_{{\rm rest}_{s-1}}}\varphi_{{\rm rest}_{s-1},[K_s]}$ has full rank and thus no components $y_{{\rm rest}_s}$ remain. At this point, i.e., after the $(s-1)$-th step, we already have
		\begin{gather*}
				y_{1,[0,\kappa_1-1]} =\varphi_{1,[0,\kappa_1-1]}(x),\\
				y_{1,[\kappa_1]} =v_{1},\\
				y_{2,[0,\kappa_2-1]} =\varphi_{2,[0,\kappa_2-1]}(x,v_{1},v_{1,[1]},\dots),\\
				y_{2,[\kappa_2]} =v_{2},\\
				 \vdotswithin{=}\\
				y_{s-1,[0,\kappa_{s-1}-1]} =\varphi_{s-1,[0,\kappa_{s-1}-1]}(x,v_{1},v_{1,[1]},\dots,v_{s-2},v_{s-2,[1]},\dots),\\
				y_{s-1,[\kappa_{s-1}]} =v_{s-1},\\
			y_{{\rm rest}_{s-1},[0,K_{{\rm rest}_{s-1}}-1]} =\varphi_{{\rm rest}_{s-1},[0,K_{{\rm rest}_{s-1}}-1]}(x,v_{1},v_{1,[1]},\dots,v_{s-2},v_{s-2,[1]},\dots),\\
				y_{{\rm rest}_{s-1},[K_{{\rm rest}_{s-1}}]} =\varphi_{{\rm rest}_{s-1},[K_{{\rm rest}_{s-1}}]}(x,v_{1},v_{1,[1]},\dots,v_{s-2},v_{s-2,[1]},\dots,v_{s-1}).
		\end{gather*}
		The multi-index $K_s=\bigl(k_s^1,\dots,k_s^{m-m_1-\dots-m_{s-1}}\bigr)$ is again defined such that
		\begin{gather*}
				\Lie_{f_u}^{k_s^j-1}\varphi_{{\rm rest}_{s-1}}^j=\varphi_{{\rm rest}_{s-1},[k_s^j-1]}^j(x,v_{1},v_{1,[1]},\dots,v_{s-1},v_{s-1,[1]},\dots),\\
				\Lie_{f_u}^{k_s^j}\varphi_{{\rm rest}_{s-1}}^j=\varphi_{{\rm rest}_{s-1},[k_s^j]}^j(x,v_{1},v_{1,[1]},\dots,v_{s-1},v_{s-1,[1]},\dots,u_{{\rm rest}_{s-1}}),
		\end{gather*}
		where by assumption we now have a regular Jacobian matrix $\partial_{u_{{\rm rest}_{s-1}}}\varphi_{{\rm rest}_{s-1},[K_s]}$ (and thus $\varphi_s=\varphi_{{\rm rest}_{s-1}}$ and $\kappa_s=K_s$). In conclusion, in the new coordinates successively constructed by this procedure, the flat output and its time derivatives up to the orders $\kappa_i$ are given by
		\begin{gather}
				y_{1,[0,\kappa_1-1]} =\varphi_{1,[0,\kappa_1-1]}(x),\nonumber\\
				y_{1,[\kappa_1]} =v_{1},\nonumber\\
				y_{2,[0,\kappa_2-1]} =\varphi_{2,[0,\kappa_2-1]}(x,v_{1},v_{1,[1]},\dots),\nonumber\\
				y_{2,[\kappa_2]} =v_{2},\nonumber\\
				 \vdotswithin{=}\nonumber\\					
y_{s-1,[0,\kappa_{s-1}-1]} =\varphi_{s-1,[0,\kappa_{s-1}-1]}(x,v_{1},v_{1,[1]},\dots,v_{s-2},v_{s-2,[1]},\dots),\nonumber\\
				y_{s-1,[\kappa_{s-1}]} =v_{s-1},\nonumber\\
				y_{s,[0,\kappa_s-1]} =\varphi_{s,[0,\kappa_s-1]}(x,v_{1},v_{1,[1]},\dots,v_{s-1},v_{s-1,[1]},\dots),\nonumber\\
				y_{s,[\kappa_s]} =\varphi_{s,[\kappa_s]}(x,v_{1},v_{1,[1]},\dots,v_{s-1},v_{s-1,[1]},\dots,u_{{\rm rest}_{s-1}})\,(=v_{s}).\label{eq:ys}
		\end{gather}
		\begin{Remark}\label{rem:procedureLinearizingOutput}
			If the procedure is applied to a flat output which is also a linearizing output in the sense of static feedback linearization, i.e., an output with a vector relative degree of $n$, see, e.g., \cite{Isidori:1995}, then due to the regularity of $\partial_u\varphi_{[R]}(x,u)$ the first step is already the last step and we have $\kappa_1=K_1=R$.
		\end{Remark}
		\begin{Remark}
			In principle, the above proposed procedure is similar to an application of the dynamic extension algorithm, which was introduced for affine-input systems in \cite{NijmeijerRespondek:1988} and is used for solving the dynamic input-output decoupling problem. A version of this algorithm for the general nonlinear case can be found in \cite{Respondek:1990}. As the extension algorithm, also our procedure essentially consists of successively replacing components of the input by certain derivatives of the components of the output. However, there are some important differences: The dynamic extension algorithm may be applied to any system with output. The outcome of the extension algorithm is then the number of decouplable input-output channels and a dynamic feedback which achieves the decoupling. In our procedure, we restrict to flat systems with an $(x,u)$-flat output. For any $(x,u)$-flat output, the outcome of our procedure are adapted coordinates for~$\XxUlu$ such that a feasible multi-index $\kappa$ can be read off. Additionally, the procedure reveals that the time derivatives of the flat output up to the orders $\kappa-1$ depend on the state~$x$ as well as higher-order time derivatives in a special triangular way. The latter will be crucial for the tracking control design in Section \ref{sec:control}.
		\end{Remark}
		\begin{Remark}
			It can be shown that one possible choice for a multi-index $A\leq R$ with $|A|=n$ such that the differentials $\D x,\D\varphi_{[A]},\D\varphi_{[A+1]},\dots,\D\varphi_{[R-1]}$ are linearly independent (i.e., a multi-index $\kappa$ in the above introduced notation) is given by the so-called structure at infinity of this $(x,u)$-flat output (after eventually permuting the components of the flat output), see, e.g., \cite{DiBenedettoGrizzleMoog:1989,GrizzleDiBenedettoMoog:1987,Martin:1993-2,Moog:1988,Nijmeijer:1986,NijmeijerSchumacher:1985,Respondek:1990} for a definition of the structure at infinity, methods for computing it, and its applications in problems like system inversion and the dynamic decoupling problem. In particular in \cite{Martin:1993-2}, an intrinsic approach to the dynamic input-output decoupling problem utilizing the structure at infinity is presented, which is applicable even to systems for which the extension algorithm does not succeed. Below, we show that our procedure applied to an~$(x,u)$-flat output always yields a multi-index $\kappa$ with the desired properties. Nevertheless, the fact that the structure at infinity represents such a multi-index $\kappa$ already implies, in combination with Theorem \ref{thm:linerizing_feedback}, that $(x,u)$-flat systems can in principle be exactly feedback linearized by a~quasi-static feedback of the form \eqref{eq:intro_quasi_static_feedback_x}. However, the structure at infinity is only one possible choice to meet the conditions of Theorem \ref{thm:linerizing_feedback}. Due to its degrees of freedom, the above proposed procedure allows to find not only the structure at infinity but also other suitable minimal multi-indices $\kappa$. This is particularly important for two reasons: First, different choices for $\kappa$ may lead to linearizing feedbacks with differently located singularities, which can be crucial in practical applications. Second, the multi-index $\kappa$ determines the orders of the individual tracking error systems of a subsequently designed tracking control (see Section \ref{sec:control}).
		\end{Remark}
		\begin{Lemma}
			For every $(x,u)$-flat output $y=\varphi(x,u)$ of the system \eqref{eq:intro_nlsys}, the stated procedure terminates after at most $m$ steps.
		\end{Lemma}
		\begin{proof}
			The existence of a multi-index $K_i=\bigl(k_i^1,\dots,k_i^{m-m_1-\dots-m_{i-1}}\bigr)$ such that $\Lie_{f_u}^{k_i^j}\varphi_{{\rm rest}_{i-1}}^j$ explicitly depends on at least one of the inputs belonging to $u_{{\rm rest}_{i-1}}$ in each step (set $\varphi_{{\rm rest}_0}=\varphi$ and~${u_{{\rm rest}_0}=u}$) is a direct consequence of the linear independence of the differentials $\D\varphi,\D\Lie_{f_u}\varphi,\allowbreak\dots,\smash{\D\Lie_{f_u}^\beta\varphi}$ of time derivatives of a flat output for arbitrary $\beta$. Indeed, assume that \smash{$\Lie_{f_u}^l\varphi_{{\rm rest}_{i-1}}^j$} for arbitrarily large $l$ only depends on $x$ and $v_1,v_{1,[1]},\dots,v_{i-1},v_{i-1,[1]},\dots$ but not on $u_{{\rm rest}_{i-1}}$. Then because of
\begin{align*}
\D\Lie_{f_u}^{l}\varphi_{{\rm rest}_{i-1}}^j&\in\Span{\D x,\D v_1,\D v_{1,[1]},\dots,\D v_{i-1},\D v_{i-1,[1]},\dots}\\
&\subset\Span{\D\varphi_{[0,R-1]},\D\varphi_1,\D\varphi_{1,[1]},\dots,\D\varphi_{i-1},\D\varphi_{i-1,[1]},\dots},
\end{align*}
 the differentials \smash{$\D\Lie_{f_u}^{l}\varphi_{{\rm rest}_{i-1}}^j$} could be expressed as a linear combination of the differentials $\D\varphi_{[0,R-1]},\D\varphi_1,\D\varphi_{1,[1]},\dots,\D\varphi_{i-1},\D\varphi_{i-1,[1]},\dots$. However, for $l\geq r$, with $r$ being the order of the highest time derivative of \smash{$y_{{\rm rest}_{i-1}}^j$} needed in the flat parameterization \eqref{eq:parameterization}, this would be a~contradiction to the linear independence of the differentials \smash{$\D\varphi,\D\Lie_{f_u}\varphi,\dots,\D\Lie_{f_u}^\beta\varphi$} for arbitrary~$\beta$. Consequently we have $m_i=\Rank{\partial_{u_{{\rm rest}_{i-1}}}\varphi_{{\rm rest}_{i-1},[K_i]}}\geq 1$ in every step, and since there are only~$m$ inputs the procedure terminates after at most $m$ steps.
		\end{proof}

		\begin{Theorem}\label{thm:construction}
			For every $(x,u)$-flat output $y=\varphi(x,u)$ of system \eqref{eq:intro_nlsys}, the multi-index $\kappa=(\kappa_1,\dots,\kappa_s)$ formed by the multi-indices $\kappa_i$ constructed in the above procedure satisfies $|\kappa|=n$ and~${\kappa\leq R}$. Furthermore, the differentials $\D x,\D\varphi_{[\kappa]},\D\varphi_{[\kappa+1]},\dots,\D\varphi_{[R-1]}$ are linearly independent.
		\end{Theorem}
		\begin{proof}
			The existence of the flat parameterization $x=F_x(\varphi_{[0,R-1]})$ implies
\[
					\Span{\D x}\subset\Span{\D\varphi_{[0,R-1]}}.
\]
			In other words, there exist exactly $n$ independent linear combinations of the differentials of the flat output and its time derivatives which are contained in $\Span{\D x}$. Now consider the expressions for the flat output and its time derivatives \eqref{eq:ys} in the coordinates successively constructed during the procedure. Because of the independence of the time derivatives of a flat output, it is possible to construct exactly $|\kappa|$ independent linear combinations of the differentials of the functions $\varphi_{[0,\kappa-1]}$ in \eqref{eq:ys} and the differentials
			\begin{align*}
					\D v_{1},\,\D v_{1,[1]},\,\D v_{1,[2]},\,\dots,\,\D v_{2},\,\D v_{2,[1]},\,\D v_{2,[2]},\,\dots,\,\D v_{s-1},\,\D v_{s-1,[1]},\,\D v_{s-1,[2]},\,\dots,
			\end{align*}
			which are contained in $\Span{\D x}$ ($|\kappa|<n$ would contradict with $y$ being a flat output, ${|\kappa|>n}$ would contradict with the linear independence of the differentials of components of time derivatives of a flat output up to arbitrary orders). Consequently, $|\kappa|=n$ follows.
			
			With $|\kappa|=n$, the property $\kappa\leq R$ can be shown by contradiction as follows. Assume that~${\kappa^j>r^j}$ for some $j\in\{1,\dots,m\}$. Then, at least one of the functions $\varphi_{[0,\kappa-1]}$ in \eqref{eq:ys} does not belong to $\varphi_{[0,R-1]}$, and consequently $\Span{\D x}\subset\Span{\D\varphi_{[0,R-1]}}$ could not hold.
			
			The linear independence of the differentials $\D x,\D\varphi_{[\kappa]},\D\varphi_{[\kappa+1]},\dots,\D\varphi_{[R-1]}$ can easily be verified in the constructed coordinates, where they are given by $\D x,\D v,\D v_{[1]},\dots,\D v_{[R-\kappa-1]}$.
		\end{proof}

		\begin{Corollary}\label{cor:procedure_dependencies}
			The time derivatives of $v$ which explicitly occur in the functions $\varphi_{i,[0,\kappa_i-1]}$, $i=1,\dots,s$ in \eqref{eq:ys} belong to the time derivatives of the $(x,u)$-flat output $y=\varphi(x,u)$ up to the orders $(R-1)$.
		\end{Corollary}
		\begin{proof}
			Assume for contradiction that there is some function among $\varphi_{i,[0,\kappa_i-1]}$, $i=1,\dots,s$ of~\eqref{eq:ys} which explicitly depends on a time derivative of $v$ (recall that $v$ are just certain time derivatives of the $(x,u)$-flat output $\varphi$) that does not belong to $\varphi_{[0,R-1]}$. Then, because of~${|\kappa|=\Dim{x}=n}$, $\Span{\D x}\subset\Span{\D\varphi_{[0,R-1]}}$ could not hold.
		\end{proof}

		According to Theorem \ref{thm:construction}, for every $(x,u)$-flat output of system \eqref{eq:intro_nlsys}, there exists a multi-index $\kappa\leq R$ with $|\kappa|=n$ such that the differentials $\D x,\D\varphi_{[\kappa]},\D\varphi_{[\kappa+1]},\dots,\D\varphi_{[R-1]}$ are linearly independent. Therefore, the conditions of Theorem \ref{thm:linerizing_feedback} are met with $A=\kappa$, and consequently the time derivatives $y_{[\kappa]}=\varphi_{[\kappa]}(x,u,u_{[1]},\dots)$ can be introduced as a new input by means of a~quasi-static feedback of the state $x$. We thus have the following corollary.
		\begin{Corollary}\label{cor:feedback_xu}
			Every $(x,u)$-flat system \eqref{eq:intro_nlsys} can be exactly feedback linearized with respect to every $(x,u)$-flat output $y=\varphi(x,u)$ by a quasi-static feedback of the form
			\begin{align}\label{eq:feedbackXU}
				\begin{aligned}
					u&=F_u\circ\Psi(x,v_{[0,R-\kappa]}).
				\end{aligned}
			\end{align}
			The input-output behavior of the closed-loop system is given by $y_{[\kappa]}=v$.
		\end{Corollary}
		\begin{Remark}\label{rem:feedback_xu}
			Note that the linearizing feedback \eqref{eq:feedbackXU} follows directly from the procedure for constructing the multi-index $\kappa$ described in this subsection. Indeed, in the coordinates successively introduced during the procedure, the flat output and its derivatives read as \eqref{eq:ys}. Feedback \eqref{eq:feedbackXU} is thus simply obtained by substituting the expressions for $y_{[0,\kappa-1]}$ from \eqref{eq:ys} and $y_{[\kappa,R]}=v_{[0,R-\kappa]}$ into the flat parameterization $F_u$ of the control input.
		\end{Remark}
		\begin{Remark}
			The main statement of Corollary \ref{cor:feedback_xu} can also be derived by showing that every $(x,u)$-flat output $y=\varphi(x,u)$ meets the condition of \cite[Theorem 2.2.2]{DelaleauPereiradaSilva:1998-2}. The motivation for providing an alternative derivation of this result within our geometric framework is that, as a~byproduct, we obtain coordinates which reveal the triangular structure \eqref{eq:ys}. The construction of these coordinates is crucial for the design of tracking control laws of the desired form~\smash{$u=\alpha\bigl(x,y^d_{[0,R]}\bigr)$}, which only depend on the systems' original state $x$ as well as the reference trajectory $y^d(t)$ and its derivatives.
		\end{Remark}	
		\section{Examples}\label{sec:examples}
		In this section, we demonstrate the procedure of Section \ref{subsec:exact_lin_xu} for the construction of linearizing quasi-static feedbacks of classical states by two $(x,u)$-flat examples. The first one is an academic example whereas the second one is a practical example. Expressions involved in the practical example are rather extensive and cannot be given in much detail. For this reason, the procedure is also illustrated by means of the simpler academic example. Furthermore, for the practical example, the proposed procedure terminates in two steps, whereas three steps are needed for the academic example.
		
	\subsection{Academic example}\label{subsec:academic}
		Consider the four-input system\footnote{This system has been constructed on basis of a three-input system which is considered in \cite{Kolar:2017} and originates from Philippe Martin.}
		\begin{alignat}{3}
				& \dot{x}^1 =u^1,\qquad && \dot{x}^6 =x^7\bigl(u^1u^3-u^2-1\bigr)+u^1x^4\bigl(u^1+x^4\bigr)-x^8u^1,&\nonumber\\
				& \dot{x}^2 =x^9,\qquad && \dot{x}^7 =x^4+u^1,&\nonumber\\
				& \dot{x}^3 =u^2-u^1u^3,\qquad &&\dot{x}^8 =x^4x^7u^1-x^6,&\nonumber\\
				& \dot{x}^4 =u^3,\qquad &&\dot{x}^9 =x^{10}+u^2+u^3,&\nonumber\\
				& \dot{x}^5 =x^3+x^4u^1, \qquad && \dot{x}^{10} =u^4 & \label{eq:academic}
		\end{alignat}
		with the $x$-flat output $y=\bigl(x^1,x^2,x^5,x^8\bigr)$. Since $x$-flat outputs are contained in the broader class of $(x,u)$-flat outputs, the results of Section~\ref{subsec:exact_lin_xu} can of course be applied to this $x$-flat output, but for demonstration purposes, let us use the $(x,u)$-flat output $y=\bigl(x^1,x^2,x^5,x^8+u^1\bigr)$ instead.
		
		The multi-index $R$ containing the highest orders of the derivatives of the components of the flat output in \eqref{eq:parameterization} is given by $R=(6,3,5,5)$. The vector field \eqref{eq:fu} for this system reads
		\begin{align}
				f_u={}&u^1\partial_{x^1}+x^9\partial_{x^2}+\bigl(u^2-u^1u^3\bigr)\partial_{x^3}+u^3\partial_{x^4}+\bigl(x^3+x^4u^1\bigr)\partial_{x^5}\nonumber\\
				&+\bigl(x^7\bigl(u^1u^3-u^2-1\bigr)+u^1x^4\bigl(u^1+x^4\bigr)-x^8u^1\bigr)\partial_{x^6}\nonumber\\
				&+\bigl(x^4+u^1\bigr)\partial_{x^7}+\bigl(x^4x^7u^1-x^6\bigr)\partial_{x^8}+\bigl(x^{10}+u^2+u^3\bigr)\partial_{x^9}\nonumber\\
				&+u^4\partial_{x^{10}}
+\sum_{\alpha=0}^{l_u-1}\bigl(u^1_{[\alpha+1]}\partial_{u^1_{[\alpha]}}+u^2_{[\alpha+1]}\partial_{u^2_{[\alpha]}}+u^3_{[\alpha+1]}\partial_{u^3_{[\alpha]}}+u^4_{[\alpha+1]}\partial_{u^4_{[\alpha]}}\bigr).\label{eq:example1_fu}
		\end{align}
		In the following, we apply the procedure stated in Section \ref{subsec:exact_lin_xu} to obtain a multi-index $\kappa$ with the properties stated in Theorem \ref{thm:construction}. Subsequently, we derive a quasi-static feedback \eqref{eq:feedbackXU} which introduces the corresponding derivatives $v=y_{[\kappa]}$ of the flat output as new inputs, and thus exactly feedback linearizes the system.
		
		\textit{Step $1$.} Differentiating the components of the flat output along the vector field \eqref{eq:example1_fu}, we find that $K_1=(1,2,1,0)$. The corresponding $k_1^j$-th derivatives of the components of the flat output are given by
		\begin{align*}
				\varphi_{[K_1]}&=\begin{bmatrix}
					\varphi^1_{[1]}\\
					\varphi^2_{[2]}\\
					\varphi^3_{[1]}\\
					\varphi^4
				\end{bmatrix}=
				\begin{bmatrix}
					u^1\\
					x^{10}+u^2+u^3\\
					x^3+x^4u^1\\
					x^8+u^1
				\end{bmatrix},
		\end{align*}
		and we have $m_1\!=\Rank{\partial_u\varphi_{[K_1]}}\!=2$. We obviously have $\Rank{\partial_u\varphi_{1,[\kappa_1]}}\!=m_1$ with ${\varphi_1\!=\bigl(\varphi^1,\varphi^2\bigr)}$, $\kappa_1=\bigl(k_1^1,k_1^2\bigr)=(1,2)$, and accordingly $\varphi_{{\rm rest}_1}=\bigl(\varphi^3,\varphi^4\bigr)$, $K_{{\rm rest}_1}=\bigl(k_1^3,k_1^4\bigr)=(1,0)$. Applying the change of coordinates
		\begin{gather*}
				v_{1}=\varphi_{1,[\kappa_1]}=\begin{bmatrix}
					\varphi^1_{[1]}\\
					\varphi^2_{[2]}
				\end{bmatrix}=\begin{bmatrix}
					u^1\\
					x^{10}+u^2+u^3
				\end{bmatrix},\\		
				v_{1,[1]}=\varphi_{1,[\kappa_1+1]}=\begin{bmatrix}
					\varphi^1_{[2]}\\
					\varphi^2_{[3]}
				\end{bmatrix}=\begin{bmatrix}
					u^1_{[1]}\\
					u^4+u^2_{[1]}+u^3_{[1]}
				\end{bmatrix},\\
				\vdotswithin{=}
		\end{gather*}
		by which we replace $u^1$ and $u^2$ and their derivatives by $v_{1}$ and its derivatives and keep $u_{{\rm rest}_1}=\bigl(u^3,u^4\bigr)$ and its derivatives as coordinates, yields
		\begin{align*}
				\varphi_{1,[\kappa_1]}=v_{1}=\begin{bmatrix}
					v_{1}^1\\
					v_{1}^2
				\end{bmatrix},
\qquad
				\varphi_{{\rm rest}_1,[K_{{\rm rest}_1}]}=\begin{bmatrix}
					\varphi^3_{[1]}\\
					\varphi^4
				\end{bmatrix}=\begin{bmatrix}
					x^3+x^4v_{1}^1\\
					x^8+v_{1}^1
				\end{bmatrix}.
		\end{align*}

		\textit{Step $2$.} We proceed by differentiating $\varphi_{{\rm rest}_1}=\bigl(\varphi^3,\varphi^4\bigr)$ until an explicit dependence on $u_{{\rm rest}_1}=\bigl(u^3,u^4\bigr)$ occurs. We have
		\begin{gather*}
				\Lie_{f_u}^2\varphi^3=v_{1}^2-x^{10}-u^3+x^4v_{1,[1]}^1,\\
				\Lie_{f_u}\varphi^4=x^4x^7v_1^1-x^6+v_{1,[1]}^1,\\
				\Lie_{f_u}^2\varphi^4=x^7\bigl(v_{1}^2-x^{10}-u^3+x^4v_{1,[1]}^1+1\bigr)+x^8v_{1}^1+v_{1,[2]}^1,
		\end{gather*}
		and thus $K_2=(2,2)$. Since only $u^3$ occurs explicitly, we have $m_2=\Rank{\partial_{u_{{\rm rest}_1}}\varphi_{{\rm rest}_1,[K_2]}}=1$. We obviously have $\Rank{\partial_{u_{{\rm rest}_1}}\varphi_{2,[\kappa_2]}}=m_2$ with $\varphi_2=\varphi^3$, $\kappa_2=k_2^1=2$, and accordingly $\varphi_{{\rm rest}_2}=\varphi^4$, $K_{{\rm rest}_2}=k_2^2=2$. Applying the change of coordinates
		\begin{gather*}
\begin{split}
&v_{2}^1=\varphi_{2,[\kappa_2]}=\varphi^3_{[2]}=v_{1}^2-x^{10}-u^3+x^4v_{1,[1]}^1,\\
&v_{2,[1]}^1=\varphi_{2,[\kappa_2+1]}=\varphi^3_{[3]}=v_{1,[1]}^2-u^4-u^3_{[1]}+u^3v_{1,[1]}^1+x^4v_{1,[2]}^1,\\
&\vdotswithin{=}
\end{split}
		\end{gather*}
		by which we replace $u^3$ and its derivatives by $v_{2}$ and its derivatives and keep $u_{{\rm rest}_2}=u^4$ and its derivatives as coordinates, yields
		\begin{gather*}
				\varphi_{2,[\kappa_2]}=v_{2}^1,\qquad
				\varphi_{{\rm rest}_2,[K_{{\rm rest}_2}]}=\varphi^4_{[2]}=x^8v_{1}^1+x^7\bigl(v_{2}^1+1\bigr)+v_{1,[2]}^1.
		\end{gather*}

		\textit{Step $3$.} We proceed by differentiating $\varphi_{{\rm rest}_2}=\varphi^4$ until an explicit dependence on $u_{{\rm rest}_2}=u^4$ occurs. Because of
		\begin{gather*}
\Lie_{f_u}^3\varphi^4=x^4x^7\bigl(v_{1}^1\bigr)^2+v_{1}^1\bigl(v_{2}^1-x^6+1\bigr)+x^4\bigl(v_{2}^1+1\bigr)+x^8v_{1,[1]}^1+x^7v_{2,[1]}^1+v_{1,[3]}^1,\\ \Lie_{f_u}^4\varphi^4=\varphi^4_{[4]}\bigl(x^4,x^6,x^7,x^8,x^{10},v_{1}^1,v_{1,[1]}^1,v_{1,[2]}^1,v_{1,[4]}^1,v_{1}^2,v_{2}^1,v_{2,[1]}^1,v_{2,[2]}^1\bigr),\\ \Lie_{f_u}^5\varphi^4=\varphi^4_{[5]}\bigl(x^4,x^6,x^7,x^8,x^{10},v_{1}^1,v_{1,[1]}^1,v_{1,[2]}^1,v_{1,[3]}^1,v_{1,[5]}^1,v_{1}^2,v_{1,[1]}^2,v_{2}^1,v_{2,[1]}^1,v_{2,[2]}^1,v_{2,[3]}^1,u^4\bigr),
		\end{gather*}
		this is the case for its $5$-th derivative, and thus $\kappa_3=5$.
		
		In conclusion, in the constructed coordinates the time derivatives of the flat output up to the orders $\kappa=(\kappa_1,\kappa_2,\kappa_3)$ with $\kappa_1=(1,2)$, $\kappa_2=2$, $\kappa_3=5$ are given by
		\begin{gather}
				y_{1,[0,\kappa_1-1]}=\begin{bmatrix}
					\varphi^1\\
					\varphi^2\\
					\varphi^2_{[1]}
				\end{bmatrix}=\begin{bmatrix}
					x^1\\
					x^2\\
					x^9
				\end{bmatrix},
\qquad
				y_{1,[\kappa_1]}=\begin{bmatrix}
					\varphi^1_{[1]}\\
					\varphi^2_{[2]}
				\end{bmatrix}=\begin{bmatrix}
					v_{1}^1\\
					v_{1}^2
				\end{bmatrix},\nonumber\\
				y_{2,[0,\kappa_2-1]}=\begin{bmatrix}
					\varphi^3\\
					\varphi^3_{[1]}
				\end{bmatrix}=\begin{bmatrix}
					x^5\\
					x^3+x^4v_{1}^1
				\end{bmatrix},\qquad
				y_{2,[\kappa_2]}=\varphi^3_{[2]}=v_{2}^1,\nonumber\\
				y_{3,[\kappa_3-1]}=\begin{bmatrix}
					\varphi^4\\
					\varphi^4_{[1]}\\
					\varphi^4_{[2]}\\
					\varphi^4_{[3]}\\
					\varphi^4_{[4]}
				\end{bmatrix}=\begin{bmatrix}
					x^8+v_1^1\\
					x^4x^7v_{1}^1-x^6+v_{1,[1]}^1\\
					x^8v_{1}^1+x^7\bigl(v_{2}^1+1\bigr)+v_{1,[2]}^1\\
					x^4x^7\bigl(v_{1}^1\bigr)^2\!+v_{1}^1\bigl(v_{2}^1\!-x^6\!+1\bigr)\!+x^4\bigl(v_{2}^1\!+1\bigr)\!+x^8v_{1,[1]}^1\!+x^7v_{2,[1]}^1\!+v_{1,[3]}^1\\
					\varphi^4_{[4]}\bigl(x^4,x^6,x^7,x^8,x^{10},v_{1}^1,v_{1,[1]}^1,v_{1,[2]}^1,v_{1,[4]}^1,v_{1}^2,v_{2}^1,v_{2,[1]}^1,v_{2,[2]}^1\bigr)
				\end{bmatrix},\nonumber\\
				y_{3,\kappa_3}=
\varphi^4_{[5]}\bigl(x^4,x^6,x^7,x^8,x^{10},v_{1}^1,v_{1,[1]}^1,v_{1,[2]}^1,v_{1,[3]}^1,v_{1,[5]}^1,v_{1}^2,v_{1,[1]}^2,v_{2}^1,v_{2,[1]}^1,v_{2,[2]}^1,v_{2,[3]}^1,u^4\bigr)\nonumber\\
\phantom{y_{3,\kappa_3}}{}=v_{3}^1.\label{eq:academic_ys}
		\end{gather}
		The linearizing quasi-static feedback \eqref{eq:feedbackXU} which introduces the derivatives $v=y_{[\kappa]}$ as inputs is obtained by substituting the above expressions \eqref{eq:academic_ys} for $y_{[0,\kappa-1]}$ and $y_{[\kappa,R]}=v_{[0,R-\kappa]}$ into the flat parameterization $F_u$ of the control input (not stated here explicitly), which yields
		\begin{gather}
				u^1=v_{1}^1,\qquad
				u^2=v_{2}^1-x^4v_{1,[1]}^1,\qquad
				u^3=v_{1}^2-x^{10}-v_{2}^1+x^4v_{1,[1]}^1,\nonumber\\
				u^4=\bar{F}_u^4
\bigl(x^4\!,x^6,x^7\!,x^8,x^{10}\!,v_{1}^1,v_{1,[1]}^1,v_{1,[2]}^1,v_{1,[3]}^1,v_{1,[5]}^1,v_{1}^2,v_{1,[1]}^2,v_{2}^1,v_{2,[1]}^1,v_{2,[2]}^1,v_{2,[3]}^1,v_{3}^1\bigr).\label{eq:academic_lin_feedback}
		\end{gather}
		The linearizing feedback \eqref{eq:academic_lin_feedback} is of the desired form \eqref{eq:intro_quasi_static_feedback_x}, i.e., it is a quasi-static feedback of the state $x$ of the system.

		\subsection{3D gantry crane}\label{subsec:gantry_crane}

		\begin{figure}[t]
			\centering
			\def\svgwidth{355pt}
			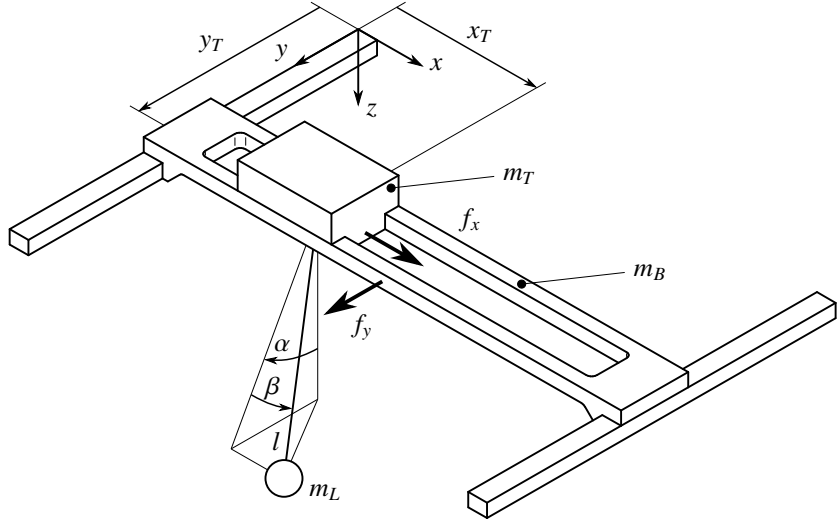
			\caption{Schematic diagram of a 3D gantry crane.}
			\label{fig:gantryCrane3D}
		\end{figure}

		Consider a gantry crane as in Figure \ref{fig:gantryCrane3D}. The trolley position is denoted by $x_T$ and $y_T$. The length of the rope is given by $l=r\phi$, where $r$ is the radius of the rope drum and $\phi$ its rotation angle. The swing angle of the rope projected onto the $yz$-plane is denoted by $\alpha$, and $\beta$ is the angle between the rope and its image under the projection onto the $yz$-plane. The load position thus follows as $x_L=x_T+r\phi\sin(\beta)$, $y_L=y_T+r\phi\sin(\alpha)\cos(\beta)$, $z_L=r\phi\cos(\alpha)\cos(\beta)$. The masses of the load, the trolley and the bridge are denoted by $m_L$, $m_T$ and $m_B$, respectively. The moment of inertia of the rope drum is denoted by $J$, and $g$ is the gravitational acceleration. The control inputs are the forces $f_x=u^1$ and $f_y=u^2$ and the torque $u^3$ acting on the rope drum. Under the assumption that the rope is always under tension, the crane can be considered as a rigid multi-body system. Its equations of motion can be derived by means of the Lagrange formalism, see, e.g., \cite{SpongVidyasagar:1989}, and follow as
		\begin{gather*}
(m_L+m_T)\ddot{x}_T+m_Lr\bigl(\sin(\beta)\ddot{\phi}+\phi\cos(\beta)\ddot{\beta}\bigr)+m_Lr\dot{\beta}\bigl(2\cos(\beta)\dot{\phi} -\phi\sin(\beta)\dot{\beta}\bigr)=u^1,\\[1mm]
(m_L+m_T+m_B)\ddot{y}_T+m_Lr\sin(\alpha)\cos(\beta)\ddot{\phi}+m_Lr\phi\bigl(\cos(\alpha)\cos(\beta)\ddot{\alpha}-\sin(\alpha)\sin(\beta)\ddot{\beta}\bigr)\\
			\qquad{}-	m_Lr\bigl(\sin(\alpha)\bigl(\phi\cos(\beta)\bigl(\dot{\alpha}^2
+\dot{\beta}^2\bigr)+2\sin(\beta)\dot{\beta}\dot{\phi}\bigr)\\
\qquad{}+2\cos(\alpha)\dot{\alpha}\bigl(\phi\sin(\beta)\dot{\beta}-\cos(\beta)\dot{\phi}\bigr)\bigr)=u^2,\\[1mm]
				m_Lr\sin(\beta)\ddot{x}_T+m_Lr\sin(\alpha)\cos(\beta)\ddot{y}_T+\bigl(J+m_Lr^2\bigr)\ddot{\phi}\\
			\qquad{}-	m_Lr\bigl(r\phi\bigl(\dot{\beta}^2+\bigl(\cos(\beta)\dot{\alpha}\bigr)^2\bigr)+g\cos(\alpha)\cos(\beta)\bigr)=u^3,\\[1mm]
				r\phi\cos(\alpha)\cos(\beta)\ddot{y}_T+\bigl(r\phi\cos(\beta)\bigr)^2\ddot{\alpha}
+r\phi\cos(\beta)\bigl(2r\dot{\alpha}\bigl(\cos(\beta)\dot{\phi}-\phi\sin(\beta)\dot{\beta}\bigr)+g\sin(\alpha)\bigr)=0,\\[1mm]
				r\phi\cos(\beta)\ddot{x}_T-r\phi\sin(\alpha)\sin(\beta)\ddot{y}_T+(r\phi)^2\ddot{\beta}\\
			\qquad{}+	r\phi\bigl(2r\dot{\beta}\dot{\phi}+r\phi\sin(\beta)\cos(\beta)\dot{\alpha}^2+\cos(\alpha)\sin(\beta)g\bigr)=0.
		\end{gather*}
		A state representation of these equations of motion reads as
		\begin{alignat}{3}
				&\dot{x}_T=v_{x_T},\qquad &&\dot{v}_{x_T}=f_{v_{x_T}}\bigl(\phi,\alpha,\beta,\omega_\alpha,\omega_\beta,u^1,u^2,u^3\bigr),&\nonumber\\
				&\dot{y}_T=v_{y_T},\qquad&& \dot{v}_{y_T}=f_{v_{y_T}}\bigl(\phi,\alpha,\beta,\omega_\alpha,\omega_\beta,u^1,u^2,u^3\bigr),&\nonumber\\
				&\dot{\phi}=\omega_\phi,\qquad&& \dot{\omega}_\phi=f_{\omega_\phi}\bigl(\phi,\alpha,\beta,\omega_\alpha,\omega_\beta,u^1,u^2,u^3\bigr),&\nonumber\\
				&\dot{\alpha}=\omega_\alpha,\qquad && \dot{\omega}_\alpha=f_{\omega_\alpha}\bigl(\phi,\alpha,\beta,\omega_\phi,\omega_\alpha,\omega_\beta,u^1,u^2,u^3\bigr),&\nonumber\\
			&	\dot{\beta}=\omega_\beta,\qquad && \dot{\omega}_\beta=f_{\omega_\beta}\bigl(\phi,\alpha,\beta,\omega_\phi,\omega_\alpha,\omega_\beta,u^1,u^2,u^3\bigr).&\label{eq:crane}
		\end{alignat}
		The load position given by
\[y=(x_T+r\phi\sin(\beta),y_T+r\phi\sin(\alpha)\cos(\beta),r\phi\cos(\alpha)\cos(\beta))\]
 forms an~$(x,u)$-flat (actually $x$-flat) output of the system. The multi-index $R$ containing the highest orders of the derivatives of the components of the flat output in \eqref{eq:parameterization} is given by $R=(4,4,4,4)$. In the following, we again apply the procedure described in Section \ref{subsec:exact_lin_xu} to obtain a~multi-index~$\kappa$ with the properties stated in Theorem \ref{thm:construction}. Subsequently, we again derive a~quasi-static feedback~\eqref{eq:intro_quasi_static_feedback_x} which introduces the corresponding derivatives $v=y_\kappa$ of the flat output as new inputs.
		
		\textit{Step $1$.} Differentiating the components of the flat output along the corresponding vector field~\eqref{eq:fu}, we find that $K_1=(2,2,2)$ with the corresponding time derivatives
		\begin{align*}
				\varphi_{[K_1]}=\begin{bmatrix}
					\varphi^1_{[2]}\bigl(\phi,\alpha,\beta,\omega_\alpha,\omega_\beta,u^1,u^2,u^3\bigr)\\
					\varphi^2_{[2]}\bigl(\phi,\alpha,\beta,\omega_\alpha,\omega_\beta,u^1,u^2,u^3\bigr)\\
					\varphi^3_{[2]}\bigl(\phi,\alpha,\beta,\omega_\alpha,\omega_\beta,u^1,u^2,u^3\bigr)
				\end{bmatrix}
		\end{align*}
		and $\Rank{\partial_u\varphi_{[K_1]}}\!=1$. We can choose, e.g., $\varphi_1\!=\varphi^3$ with $\kappa_1\!=k_1^3\!=2$ and thus ${\varphi_{{\rm rest}_1}\!=\bigl(\varphi^1,\varphi^2\bigr)}$, and apply the change of coordinates
		\begin{gather*}
				v_{1}^1=\varphi_{1,[\kappa_1]}=\varphi^3_{[2]}\bigl(\phi,\alpha,\beta,\omega_\alpha,\omega_\beta,u^1,u^2,u^3\bigr),\\
				v_{1,[1]}^1=\varphi_{1,[\kappa_1+1]}=\varphi^3_{[3]}\bigl(\phi,\alpha,\beta,\omega_\alpha,\omega_\beta,\omega_\phi,u^1,u^2,u^3,u^1_{[1]},u^2_{[1]},u^3_{[1]}\bigr),\\
				\vdotswithin{=}
		\end{gather*}
		by which we replace the input $u^3$ and its derivatives by $v_{1}$ and its derivatives, and keep ${u_{{\rm rest}_1}=\bigl(u^1,u^2\bigr)}$ and its derivatives as coordinates. This results in
		\begin{gather*}
				\varphi_{1,[\kappa_1]}=v_{1}^1,\qquad
				\varphi_{{\rm rest}_1,[K_{{\rm rest}_1}]}=\begin{bmatrix}
					\varphi^1_{[2]}\bigl(\alpha,\beta,v_{1}^1\bigr)\\
					\varphi^2_{[2]}\bigl(\alpha,v_{1}^1\bigr)
				\end{bmatrix}.
		\end{gather*}
		Note that the choice $\varphi_1=\varphi^3$ is not unique but the most practical one, since the other two possible choices $\varphi_1=\varphi^1$ and $\varphi_1=\varphi^2$, would lead to feedback laws with singularities for $\alpha=0$ or $\beta=0$.
		
		\textit{Step $2$.} We proceed by differentiating $\varphi_{{\rm rest}_1}=\bigl(\varphi^1,\varphi^2\bigr)$ until an explicit dependence on $u_{{\rm rest}_1}=\bigl(u^1,u^2\bigr)$ occurs. This is the case for the $4$-th derivatives
		\begin{gather*}
				\Lie_{f_u}^4\varphi^1=\varphi^1_{[4]}\bigl(\phi,\alpha,\beta,\omega_\phi,\omega_\alpha,\omega_\beta,v_{1}^1,v_{1,[1]}^1,v_{1,[2]}^1,u^1\bigr),\\
				\Lie_{f_u}^4\varphi^2=\varphi^2_{[4]}\bigl(\phi,\alpha,\beta,\omega_\phi,\omega_\alpha,\omega_\beta,v_{1}^1,v_{1,[1]}^1,v_{1,[2]}^1,u^2\bigr),
		\end{gather*}
		i.e., $K_2=(4,4)$. Because of $m_2=\Rank{\partial_{u_{{\rm rest}_1}}\varphi_{{\rm rest}_1,[K_2]}}=2$, the procedure terminates at this point and we have $\kappa_2=K_2$.
		
		In conclusion, in the constructed coordinates the time derivatives of the flat output up to the orders $\kappa=(\kappa_1,\kappa_2)$ with $\kappa_1=2$, $\kappa_2=(4,4)$ are of the form
		\begin{gather}
				y_{1,[0,\kappa_1-1]}=\begin{bmatrix}
					\varphi^3(\phi,\alpha,\beta)\\
					\varphi^3_{[1]}(\phi,\alpha,\beta,\omega_\phi,\omega_\alpha,\omega_\beta)
				\end{bmatrix},\qquad
				y_{1,[\kappa_1]}=\varphi^3_{[2]}=v_{1}^1,\nonumber\\
				y_{2,[0,\kappa_2-1]}=\begin{bmatrix}
					\varphi^1(x_T,\phi,\beta)\\
					\varphi^1_{[1]}(v_{x_T},\phi,\beta,\omega_\phi,\omega_\beta)\\
					\varphi^1_{[2]}\bigl(\alpha,\beta,v_{1}^1\bigr)\\
					\varphi^1_{[3]}\bigl(\alpha,\beta,\omega_\alpha,\omega_\beta,v_{1}^1,v_{1,[1]}^1\bigr)\\
					\varphi^2(y_T,\phi,\alpha,\beta)\\
					\varphi^2_{[1]}(\phi,\alpha,\beta,v_{y_T},\omega_\phi,\omega_\alpha,\omega_\beta)\\
					\varphi^2_{[2]}\bigl(\alpha,v_{1}^1\bigr)\\
					\varphi^2_{[3]}\bigl(\alpha,\omega_\alpha,v_{1}^1,v_{1,[1]}^1\bigr)
				\end{bmatrix},\nonumber\\
				y_{2,\kappa_2}=\begin{bmatrix}
					\varphi^1_{[4]}\bigl(\phi,\alpha,\beta,\omega_\phi,\omega_\alpha,\omega_\beta,v_{1}^1,v_{1,[1]}^1,v_{1,[2]}^1,u^1\bigr)\\
					\varphi^2_{[4]}\bigl(\phi,\alpha,\beta,\omega_\phi,\omega_\alpha,\omega_\beta,v_{1}^1,v_{1,[1]}^1,v_{1,[2]}^1,u^2\bigr)
				\end{bmatrix}=\begin{bmatrix}
					v_{2}^1\\
					v_{2}^2
				\end{bmatrix}.\label{eq:crane_y_e}
		\end{gather}
		The quasi-static feedback \eqref{eq:feedbackXU} which introduces the derivatives $v=y_{[\kappa]}$ as inputs is obtained by substituting the above expressions \eqref{eq:crane_y_e} for $y_{[0,\kappa-1]}$ and $y_{[\kappa,R]}=v_{[0,R-\kappa]}$ into the flat parameterization $F_u$ of the control input (not stated here explicitly), which yields a linearizing feedback of the form
		\begin{gather}
				u^1=\bar{F}_u^1\bigl(\phi,\alpha,\beta,\omega_\phi,\omega_\alpha,\omega_\beta,v_{1}^1,v_{1,[1]}^1,v_{1,[2]}^1,v_{2}^1\bigr),\nonumber\\
				u^2=\bar{F}_u^2\bigl(\phi,\alpha,\beta,\omega_\phi,\omega_\alpha,\omega_\beta,v_{1}^1,v_{1,[1]}^1,v_{1,[2]}^1,v_{2}^2\bigr),\nonumber\\
				u^3=\bar{F}_u^3\bigl(\phi,\alpha,\beta,\omega_\phi,\omega_\alpha,\omega_\beta,v_{1}^1,v_{1,[1]}^1,v_{1,[2]}^1,v_{2}^1,v_{2}^2\bigr).\label{eq:crane_lin_feedback}
		\end{gather}
		This feedback is again of the desired form \eqref{eq:intro_quasi_static_feedback_x}.

\section[Tracking control design for (x,u)-flat systems]{Tracking control design for $\boldsymbol{(x,u)}$-flat systems}\label{sec:control}
	Trajectory tracking is a fundamental control engineering problem. Given a reference trajectory~$y^d(t)$ for the flat output $y$, the goal is to design a feedback control law which assures that the tracking error $e=y-y^d$ asymptotically decays to zero. The exact feedback linearization of a flat system, which we have discussed in Section \ref{sec:exact_lin}, is only the first step in the design of a~flatness-based tracking control. The input-output behavior of the exactly feedback linearized system is given by $m$ integrator chains, and hence the second step consists in the construction of an additional feedback such that the tracking error with respect to a reference trajectory~$y^d(t)$ is stabilized asymptotically. If the time derivatives of the flat output which form these integrator chains are available from measurements or estimates, then designing such a stabilizing feedback is of course a~straightforward task. Furthermore, if the exact feedback linearization has been performed by a~quasi-static feedback \eqref{eq:intro_quasistatic_feedback} of a generalized Brunovsk\'y state, then the time derivatives of the flat output which are needed for a stabilizing feedback are obviously given by the same generalized Brunovsk\'y state. This standard approach is discussed in detail, e.g., in~\cite{DelaleauRudolph:1998}. However, since the time derivatives of a flat output which form a generalized Brunovsk\'y state are often much harder to measure than the components of a classical state, the aim of the present paper is to provide a complete solution for the tracking control problem which requires only measurements of the latter. In the following, we show that for $(x,u)$-flat outputs \eqref{eq:xu_flat_output_procedure} it is indeed always possible to achieve a linear, decoupled and asymptotically stable tracking error dynamics with arbitrary eigenvalues by a control law of the form \smash{$u=\alpha\bigl(x,y^d_{[0,R]}(t)\bigr)$}, which only depends on the classical state $x$ and the reference trajectory $y^d(t)$. As a basis for the tracking control design we use the exact feedback linearization according to Corollary \ref{cor:feedback_xu}, which is achieved by a quasi-static feedback
	\begin{align}\label{eq:quasi_static_feedback_xu}
		\begin{aligned}
			u&=F_u\circ\Psi(x,v_{[0,R-\kappa]})=\bar{F}_u(x,v_{[0,R-\kappa]})
		\end{aligned}
	\end{align}
	of $x$. Recall that the multi-index $\kappa$ is obtained by applying the procedure stated in Section~\ref{subsec:exact_lin_xu}, and that the linearizing feedback \eqref{eq:quasi_static_feedback_xu} then follows directly by substituting the expressions for~$y_{[0,\kappa-1]}$ from \eqref{eq:ys} and $y_{[\kappa,R]}=v_{[0,R-\kappa]}$ into the flat parameterization $F_u$ of the input (see also Corollary \ref{cor:feedback_xu} and Remark \ref{rem:feedback_xu}). In the following, the special structure \eqref{eq:ys} of the derivatives of the flat output in the coordinates successively introduced during the procedure will be crucial for the tracking control design.
	
	The feedback-modified system obtained by applying feedback \eqref{eq:quasi_static_feedback_xu} has a linear input-output behavior \smash{$y_{[\kappa]}=v$} between the new input $v$ and the flat output $y$. Given a sufficiently often differentiable reference trajectory $y^d(t)$, the control law
	\begin{gather}\label{eq:control_quasistatic_feedback_x_v}
v_{i}^{j_i}=y^{{j_i},d}_{i,[\kappa_i^{j_i}]}-\sum_{\beta=0}^{\kappa_i^{j_i}-1}a_i^{{j_i},\beta}\bigl(y^{j_i}_{i,[\beta]}-y^{{j_i},d}_{i,[\beta]}\bigr),
\qquad i=1,\dots,s,\quad j_i=1,\dots,m_i,
	\end{gather}
	for the new input $v$ results in the linear tracking error dynamics
	\begin{align*}
		e^{j_i}_{i,[\kappa_i^{j_i}]}+\sum_{\beta=0}^{\kappa_i^{j_i}-1}a_i^{{j_i},\beta}e^{j_i}_{i,[\beta]}=0,
	\end{align*}
	where $e_i^{j_i}=y_i^{j_i}-y_i^{{j_i},d}$. The roots of the characteristic polynomials of the tracking error dynamics can be adjusted by the coefficients \smash{$a_i^{{j_i},\beta}\in\mathbb{R}$}. The time derivatives $y_{[0,\kappa-1]}$ occurring in \eqref{eq:control_quasistatic_feedback_x_v} can be expressed in terms of $x$ and $v_{[0,R-\kappa-1]}$, see \eqref{eq:ys} and Corollary \ref{cor:procedure_dependencies}. In \eqref{eq:quasi_static_feedback_xu}, not only $v$, but also its time derivatives $v_{[1,R-\kappa]}$ occur. Differentiating \eqref{eq:control_quasistatic_feedback_x_v} with respect to time yields
	\begin{gather*}		v_{i,[\lambda]}^{j_i}=y^{{j_i},d}_{i,[\kappa_i^{j_i}+\lambda]}-\sum_{\beta=0}^{\kappa_i^{j_i}-1}a_i^{{j_i},\beta}\bigl(y^{j_i}_{i,[\beta+\lambda]}-y^{{j_i},d}_{i,[\beta+\lambda]}\bigr),
\qquad\lambda=1,\dots,r_i^{j_i}-\kappa_i^{j_i}.
	\end{gather*}
	Replacing the time derivatives $y_{[0,\kappa-1]}$ by the corresponding expressions of $x$ and $v_{[0,R-\kappa-1]}$ from~\eqref{eq:ys}, and replacing $y_{[\kappa,R-1]}$ by $v_{[0,R-\kappa-1]}$ results in the system of equations
	\begin{gather}
v_{1}^{j_1} =y^{{j_1},d}_{1,[\kappa_1^{j_1}]}-\sum_{\beta=0}^{\kappa_1^{j_1}-1}a_1^{{j_1},\beta}\bigl(\varphi^{j_1}_{1,[\beta]}(x)-y^{{j_1},d}_{1,[\beta]}\bigr),\qquad{j_1}=1,\dots,m_1,\nonumber\\
			v_{1,[1]}^{j_1} =y^{{j_1},d}_{1,[\kappa_1^{j_1}+1]}-a_1^{j_1,\kappa_1^{j_1}-1}\bigl(v_{1}^{j_1}-y^{{j_1},d}_{1,[\kappa_1^{j_1}]}\bigr)
-\sum_{\beta=0}^{\kappa_1^{j_1}-2}a_1^{{j_1},\beta}\bigl(\varphi^{j_1}_{1,[\beta+1]}(x)-y^{{j_1},d}_{1,[\beta+1]}\bigr),\nonumber\\
			 \vdotswithin{=}\nonumber\\
			v_{2}^{j_2} =y^{{j_2},d}_{2,[\kappa_2^{j_2}]}-\sum_{\beta=0}^{\kappa_2^{j_2}-1}a_2^{{j_2},\beta}
\bigl(\varphi^{j_2}_{2,[\beta]}(x,v_{1},v_{1,[1]},\dots)-y^{{j_2},d}_{2,[\beta]}\bigr),\qquad{j_2}=1,\dots,m_2,\nonumber\\
			v_{2,[1]}^{j_2} =y^{{j_2},d}_{2,[\kappa_2^{j_2}+1]}-a_2^{j_2,\kappa_2^{j_2}-1}\bigl(v_{2}^{j_2}-y^{{j_2},d}_{2,[\kappa_2^{j_2}]}\bigr)
-\sum_{\beta=0}^{\kappa_2^{j_2}-2}a_2^{{j_2},\beta}\bigl(\varphi^{j_2}_{2,[\beta+1]}(x,v_{1},v_{1,[1]},\dots)-y^{{j_2},d}_{2,[\beta+1]}\bigr),\nonumber\\
			 \vdotswithin{=}\nonumber\\
			v_{3}^{j_3} =y^{{j_3},d}_{3,[\kappa_3^{j_3}]}
-\sum_{\beta=0}^{\kappa_3^{j_3}-1}a_3^{{j_3},\beta}\bigl(\varphi^{j_3}_{3,[\beta]}(x,v_{1},v_{1,[1]},\dots,v_{2},v_{2,[1]},\dots)-y^{{j_3},d}_{3,[\beta]}\bigr),\qquad{j_3}=1,\dots,m_3,\nonumber\\
			v_{3,[1]}^{j_3} =\cdots,\nonumber\\
			 \vdotswithin{=}\nonumber\\
			v_{s}^{j_s} =y^{{j_s},d}_{s,[\kappa_s^{j_s}]}
-\sum_{\beta=0}^{\kappa_s^{j_s}-1}a_s^{{j_s},\beta}\bigl(\varphi^{j_s}_{s,[\beta]}(x,v_{1},v_{1,[1]},\dots,v_{s-1},v_{s-1,[1]},\dots)-y^{{j_s},d}_{s,[\beta]}\bigr),\nonumber\\
{j_s}=1,\dots,m_s.\label{eq:control_quasistatic_feedback_x_v_subs}
	\end{gather}
	Due to the special structure of the time derivatives $y_{[0,\kappa-1]}$ of the flat output in \eqref{eq:ys}, the time derivatives $v_{[0,R-\kappa]}$ occur in equations \eqref{eq:control_quasistatic_feedback_x_v_subs} in a triangular manner. The right-hand sides of the equations for $v_1$ do not depend on $v$ at all, the right-hand sides of the equations for $v_{1,[\lambda]}$ only depend on $v_{1,[\lambda-1]}$. The right-hand sides of the equations for $v_2$ only depend on $v_1$ and its time derivatives, and so on. The complete set of time derivatives $v_{[0,R-\kappa]}$, which is needed for practically realizing the control law \eqref{eq:control_quasistatic_feedback_x_v} for the input $v$ via the actual control input $u$, i.e., via feedback \eqref{eq:quasi_static_feedback_xu}, can thus be determined systematically by solving \eqref{eq:control_quasistatic_feedback_x_v_subs} from top to bottom for the time derivatives $v_{i,[\lambda]}^{j_i}$ as functions of $x$ and $y^d_{[0,R]}$, i.e., $v_{[0,R-\kappa]}=\phi\bigl(x,y^d_{[0,R]}\bigr)$. Substituting this solution into feedback \eqref{eq:quasi_static_feedback_xu} yields a tracking control law of the desired form
	\begin{align}\label{eq:controllaw x}
		u&=\alpha\bigl(x,y^d_{[0,R]}\bigr),
	\end{align}
	which only depends on the state $x$ and the reference trajectory.
	\begin{Remark}
		A tracking control law of the form \eqref{eq:controllaw x} has already been derived in \cite{MartinMurrayRouchon:1997} for a~particular flat system, namely the well-known PVTOL aircraft. Another derivation of such a~control law, also for the PVTOL, can be found in \cite{Rudolph:2021}. However, these derivations do not follow a systematic approach which can easily be transferred to other systems, whereas the approach presented in this paper can easily be applied to any $(x,u)$-flat system.
	\end{Remark}
		
\section{Examples continued}\label{sec:examples_continued}
	In this section, we derive tracking control laws for the two systems that were already exactly feedback linearized in Section \ref{sec:examples}.
	\subsection{Academic example}
		Consider again system \eqref{eq:academic} with the flat output $y=(y_{1},y_{2},y_{3})$ where $y_{1}=\bigl(y^1,y^2\bigr)=\bigl(x^1,x^2\bigr)$, $y_{2}=y^3=x^5$, $y_{3}=y^4=x^8+u^1$. In Section \ref{subsec:academic}, we have derived the linearizing feedback \eqref{eq:academic_lin_feedback}, which introduces $v=y_{[\kappa]}$ with $\kappa=(\kappa_1,\kappa_2,\kappa_3)$ and $\kappa_1=(1,2)$, $\kappa_2=2$, $\kappa_3=5$ as new input. The control law \eqref{eq:control_quasistatic_feedback_x_v} for $v=(v_{1},v_{2},v_{3})$ with $v_{1}=\bigl(v_{1}^1,v_{1}^2\bigr)$, $v_{2}=v_{2}^1$ and $v_{3}=v_{3}^1$ is thus given by
		\begin{align*}
			v_{1}^1={}&y_{1,[1]}^{1,d}-a_1^{1,0}\bigl(y_{1}^1-y_{1}^{1,d}\bigr),\\
			v_{1}^2={}&y_{1,[2]}^{2,d}-a_1^{2,0}\bigl(y_{1}^2-y_{1}^{2,d}\bigr)-a_1^{2,1}\bigl(y_{1,[1]}^2-y_{1,[1]}^{2,d}\bigr),\\
			v_{2}^1={}&y_{2,[2]}^{1,d}-a_2^{1,0}\bigl(y_{2}^1-y_{2}^{1,d}\bigr)-a_2^{1,1}\bigl(y_{2,[1]}^1-y_{2,[1]}^{1,d}\bigr),\\
			v_{3}^1={}&y_{3,[5]}^{1,d}-a_3^{1,0}\bigl(y_{3}^1-y_{3}^{1,d}\bigr)-a_3^{1,1}\bigl(y_{3,[1]}^1-y_{3,[1]}^{1,d}\bigr)-a_3^{1,2}\bigl(y_{3,[2]}^1-y_{3,[2]}^{1,d}\bigr)\\
			&-a_3^{1,3}\bigl(y_{3,[3]}^1-y_{3,[3]}^{1,d}\bigr)-a_3^{1,4}\bigl(y_{3,[4]}^1-y_{3,[4]}^{1,d}\bigr).
		\end{align*}
		Substituting the corresponding expressions \eqref{eq:academic_ys} for $y_{[0,\kappa-1]}$ into this control law yields
		\begin{align*}
\begin{aligned}
			v_{1}^1={}&y^{1,d}_{[1]}-a_1^{1,0}\bigl(x^1-y^{1,d}\bigr),\\
			v_{1}^2={}&y^{2,d}_{[2]}-a_1^{2,0}\bigl(x^2-y^{2,d}\bigr)-a_1^{2,1}\bigl(x^9-y^{2,d}_{[1]}\bigr),\\
			v_{2}^1={}&y^{3,d}_{[2]}-a_2^{1,0}\bigl(x^5-y^{3,d}\bigr)-a_2^{1,1}\bigl(x^3+x^4v_{1}^1-y^{3,d}_{[1]}\bigr),\\
			v_{3}^1={}&y^{4,d}_{[5]}-a_3^{1,0}\bigl(x^8+v_1^1-y^{4,d}\bigr)-a_3^{1,1}\bigl(x^4x^7v_{1}^1-x^6+v_{1,[1]}^1-y^{4,d}_{[1]}\bigr)\\
			&-a_3^{1,2}\bigl(x^8v_{1}^1+x^7\bigl(v_{2}^1+1\bigr)+v_{1,[2]}^1-y^{4,d}_{[2]}\bigr)\\
			&-a_3^{1,3}\bigl(\varphi^4_{[3]}\bigl(x^4,x^6,x^7,x^8,v_{1}^1,v_{1,[1]}^1,v_{1,[3]}^1,v_{2}^1,v_{2,[1]}^1\bigr)-y^{4,d}_{[3]}\bigr)\\
&-a_3^{1,4}\bigl(\varphi^4_{[4]}\bigl(x^4,x^6,x^7,x^8,x^{10},v_{1}^1,v_{1,[1]}^1,v_{1,[2]}^1,v_{1,[4]}^1,v_{1}^2,v_{2}^1,v_{2,[1]}^1,v_{2,[2]}^1\bigr)
-y^{4,d}_{[4]}\bigr).
\end{aligned}
		\end{align*}
		In order to express the required time derivatives $v_{[0,R-\kappa]}$ in terms of $x$ and $y^d_{[0,R]}$, we differentiate the equations for $v$ up to the orders $R-\kappa$, which yields
		\begin{gather*}
			v_{1}^1 =y^{1,d}_{[1]}-a_1^{1,0}\bigl(x^1-y^{1,d}\bigr),\\
			v_{1,[1]}^1=y^{1,d}_{[2]}-a_1^{1,0}\bigl(v_{1}^1-y^{1,d}_{[1]}\bigr),\\
			v_{1,[2]}^1 =y^{1,d}_{[3]}-a_1^{1,0}\bigl(v_{1,[1]}^1-y^{1,d}_{[2]}\bigr),\\
			v_{1,[3]}^1=y^{1,d}_{[4]}-a_1^{1,0}\bigl(v_{1,[2]}^1-y^{1,d}_{[3]}\bigr),\\
			v_{1,[4]}^1 =y^{1,d}_{[5]}-a_1^{1,0}\bigl(v_{1,[3]}^1-y^{1,d}_{[4]}\bigr),\\
			v_{1,[5]}^1=y^{1,d}_{[6]}-a_1^{1,0}\bigl(v_{1,[4]}^1-y^{1,d}_{[5]}\bigr),\\
			v_{1}^2 =y^{2,d}_{[2]}-a_1^{2,0}\bigl(x^2-y^{2,d}\bigr)-a_1^{2,1}\bigl(x^9-y^{2,d}_{[1]}\bigr),\\
			v_{1,[1]}^2 =y^{2,d}_{[3]}-a_1^{2,0}\bigl(x^9-y^{2,d}_{[1]}\bigr)-a_1^{2,1}\bigl(v_1^2-y^{2,d}_{[2]}\bigr),\\
			v_{2}^1 =y^{3,d}_{[2]}-a_2^{1,0}\bigl(x^5-y^{3,d}\bigr)-a_2^{1,1}\bigl(x^3+x^4v_{1}^1-y^{3,d}_{[1]}\bigr),\\
			v_{2,[1]}^1 =y^{3,d}_{[3]}-a_2^{1,0}\bigl(x^3+x^4v_{1}^1-y^{3,d}_{[1]}\bigr)-a_2^{1,1}\bigl(v_{2}^1-y^{3,d}_{[2]}\bigr),\\
			v_{2,[2]}^1 =y^{3,d}_{[4]}-a_2^{1,0}\bigl(v_{2}^1-y^{3,d}_{[2]}\bigr)-a_2^{1,1}\bigl(v_{2,[1]}^1-y^{3,d}_{[3]}\bigr),\\
			v_{2,[3]}^1 =y^{3,d}_{[5]}-a_2^{1,0}\bigl(v_{2,[1]}^1-y^{3,d}_{[3]}\bigr)-a_2^{1,1}\bigl(v_{2,[2]}^1-y^{3,d}_{[4]}\bigr),\\
			v_{3}^1 =y^{4,d}_{[5]}-a_3^{1,0}\bigl(x^8+v_{1}^1-y^{4,d}\bigr)-a_3^{1,1}\bigl(x^4x^7v_{1}^1-x^6+v_{1,[1]}^1-y^{4,d}_{[1]}\bigr)\\
\phantom{v_{3}^1 =}{}	- a_3^{1,2}\bigl(x^8v_{1}^1+x^7(v_{2}^1+1)+v_{1,[2]}^1-y^{4,d}_{[2]}\bigr)\\
\phantom{v_{3}^1 =}{}	- a_3^{1,3}\bigl(\varphi^4_{[3]}\bigl(x^4,x^6,x^7,x^8,v_{1}^1,v_{1,[1]}^1,v_{1,[3]}^1,v_{2}^1,v_{2,[1]}^1\bigr)-y^{4,d}_{[3]}\bigr)\\
\phantom{v_{3}^1 =}{}	- a_3^{1,4}\bigl(\varphi^4_{[4]}\bigl(x^4,x^6,x^7,x^8,x^{10},v_{1}^1,v_{1,[1]}^1,v_{1,[2]}^1,v_{1,[4]}^1,v_{1}^2,v_{2}^1,v_{2,[1]}^1,v_{2,[2]}^1\bigr)
-y^{4,d}_{[4]}\bigr).
		\end{gather*}
		This system of equations can easily be solved from top to bottom in order to obtain $v_{[0,R-\kappa]}$ in terms of $x$ and \smash{$y^d_{[0,R]}$}. Inserting the solution into the linearizing feedback \eqref{eq:academic_lin_feedback} yields a tracking control law of the form
		\begin{align*}
			u^1={}&y^{1,d}_{[1]}-a_1^{1,0}\bigl(x^1-y^{1,d}\bigr),\\
			u^2={}&y^{3,d}_{[2]}-a_2^{1,0}\bigl(x^5-y^{3,d}\bigr)-a_2^{1,1}\bigl(x^3+x^4\bigl(y^{1,d}_{[1]}-a_1^{1,0}\bigl(x^1-y^{1,d}\bigr)\bigr)-y^{3,d}_{[1]}\bigr)\\
			&-x^4\bigl(y^{1,d}_{[2]}-a_1^{1,0}\bigl(y^{1,d}_{[1]}-a_1^{1,0}\bigl(x^1-y^{1,d}\bigr)-y^{1,d}_{[1]}\bigr)\bigr),\\
			u^3={}&y^{2,d}_{[2]}-a_1^{2,0}\bigl(x^2-y^{2,d}\bigr)-a_1^{2,1}\bigl(x^9-y^{2,d}_{[1]}\bigr)-x^{10}-y^{3,d}_{[2]}+a_2^{1,0}\bigl(x^5-y^{3,d}\bigr)\\
			&+a_2^{1,1}\bigl(x^3+x^4\bigl(y^{1,d}_{[1]}-a_1^{1,0}\bigl(x^1-y^{1,d}\bigr)\bigr)-y^{3,d}_{[1]}\bigr)\\
&+x^4\bigl(y^{1,d}_{[2]}-a_1^{1,0}\bigl(y^{1,d}_{[1]}-a_1^{1,0}\bigl(x^1-y^{1,d}\bigr)-y^{1,d}_{[1]}\bigr)\bigr),\\
			u^4={}&\alpha^4\bigl(x^1,\dots,x^{10},y^{1,d}_{[0,6]},y^{2,d}_{[0,3]},y^{3,d}_{[0,5]},y^{4,d}_{[0,5]}\bigr).
		\end{align*}
	\subsection{3D gantry crane}
		Consider again system \eqref{eq:crane} with the flat output $y=(y_1,y_2)$ where $y_1=y^3=r\phi\cos(\alpha)\cos(\beta)$, $y_2=\bigl(y^1,y^2\bigr)=(x_T+r\phi\sin(\beta),y_T+r\phi\sin(\alpha)\cos(\beta))$. In Section \ref{subsec:gantry_crane}, we have derived the linearizing feedback \eqref{eq:crane_lin_feedback}, which introduces $v=y_{[\kappa]}$ with $\kappa=(\kappa_1,\kappa_2)$ and $\kappa_1=2$, $\kappa_2=(4,4)$ as new input. The control law \eqref{eq:control_quasistatic_feedback_x_v} for $v$ with $v_{1}=v_{1}^1$, $v_{2}=\bigl(v_2^1,v_2^2\bigr)$ is thus given by
		\begin{align*}
			v_1^1&=y^{1,d}_{1,[2]}-a_1^{1,0}\bigl(y^1_1-y^{1,d}_1\bigr)-a_1^{1,1}\bigl(y^1_{1,[1]}-y^{1,d}_{1,[1]}\bigr),\\	v_2^1&=y^{1,d}_{2,[4]}-a_2^{1,0}\bigl(y^1_2-y^{1,d}_2\bigr)-a_2^{1,1}\bigl(y^1_{2,[1]}-y^{1,d}_{2,[1]}\bigr)-a_2^{1,2}\bigl(y^1_{2,[2]}-y^{1,d}_{2,[2]}\bigr)-a_2^{1,3}\bigl(y^1_{2,[3]}-y^{1,d}_{2,[3]}\bigr),\\	v_2^2&=y^{2,d}_{2,[4]}-a_2^{2,0}\bigl(y^2_2-y^{2,d}_2\bigr)-a_2^{2,1}\bigl(y^2_{2,[1]}-y^{2,d}_{2,[1]}\bigr)-a_2^{2,2}\bigl(y^2_{2,[2]}-y^{2,d}_{2,[2]}\bigr)-a_2^{2,3}\bigl(y^2_{2,[3]}-y^{2,d}_{2,[3]}\bigr).
		\end{align*}
		Substituting the corresponding expressions \eqref{eq:crane_y_e} for $y_{[0,\kappa-1]}$ into this control law and differentiating the equations for $v$ up to the corresponding orders $R-\kappa$ yields
		\begin{gather*}
			v_{1}^1 =y^{3,d}_{[2]}-a_1^{1,0}\bigl(\varphi^3(\phi,\alpha,\beta)-y^{3,d}\bigr)-a_1^{1,1}\bigl(\varphi^3_{[1]}(\phi,\alpha,\beta,\omega_\phi,\omega_\alpha,\omega_\beta)-y^{3,d}_{[1]}\bigr),\\
v_{1,{[1]}}^1 =y^{3,d}_{[3]}-a_1^{1,0}\bigl(\varphi^3_{[1]}(\phi,\alpha,\beta,\omega_\phi,\omega_\alpha,\omega_\beta)-y^{3,d}_{[1]}\bigr)-a_1^{1,1}\bigl(v_{1}^1-y^{3,d}_{[2]}\bigr),\\
			v_{1,[2]}^1 =y^{3,d}_{[4]}-a_1^{1,0}\bigl(v_{1}^1-y^{3,d}_{[2]}\bigr)-a_1^{1,1}\bigl(v_{1,[1]}^1-y^{3,d}_{[3]}\bigr),\\		
	v_{2}^1 =y^{1,d}_{[4]}-a_2^{1,0}\bigl(\varphi^1(x_T,\phi,\beta)-y^{1,d}\bigr)-a_2^{1,1}\bigl(\varphi^1_{[1]}(v_{x_T},\phi,\beta,\omega_\phi,\omega_\beta)-y^{1,d}_{[1]}\bigr)\\
\phantom{v_{2}^1 =}{} -a_2^{1,2}\bigl(\varphi^1_{[2]}\bigl(\alpha,\beta,v_{1}^1\bigr)-y^{1,d}_{[2]}\bigr)-a_2^{1,3}\bigl(\varphi^1_{[3]}\bigl(\alpha,\beta,\omega_\alpha,\omega_\beta,v_{1}^1,v_{1,[1]}^1\bigr)-y^{1,d}_{[3]}\bigr),\\
			v_{2}^2 =y^{2,d}_{[4]}-a_2^{2,0}\bigl(\varphi^2(y_T,\phi,\alpha,\beta)-y^{2,d}\bigr)-a_2^{2,1}\bigl(\varphi^2_{[1]}(\phi,\alpha,\beta,v_{y_T},\omega_\phi,\omega_\alpha,\omega_\beta)-y^{2,d}_{[1]}\bigr)\\
\phantom{v_{2}^2=}{}			-a_2^{2,2}\bigl(\varphi^2_{[2]}\bigl(\alpha,v_{1}^1\bigr)-y^{2,d}_{[2]}\bigr)-a_2^{2,3}\bigl(\varphi^2_{[3]}\bigl(\alpha,\omega_\alpha,v_{1}^1,v_{1,[1]}^1\bigr)-y^{2,d}_{[3]}\bigr).
		\end{gather*}
		This system of equations can again easily be solved from top to bottom in order to obtain~$v_{[0,R-\kappa]}$ in terms of $x$ and \smash{$y^d_{[0,R]}$}. Inserting the solution into the linearizing feedback \eqref{eq:crane_lin_feedback} yields a tracking control law of the form
		\begin{align*}
			u^1&=\alpha^1\bigl(x_T,\phi,\alpha,\beta,v_{x_T},\omega_\phi,\omega_\alpha,\omega_\beta,y^{1,d}_{[0,4]},y^{3,d}_{[0,4]}\bigr),\\
			u^2&=\alpha^2\bigl(y_T,\phi,\alpha,\beta,v_{y_T},\omega_\phi,\omega_\alpha,\omega_\beta,y^{2,d}_{[0,4]},y^{3,d}_{[0,4]}\bigr),\\
			u^3&=\alpha^3\bigl(x_T,y_T,\phi,\alpha,\beta,v_{x_T},v_{y_T},\omega_\phi,\omega_\alpha,\omega_\beta,y^{1,d}_{[0,4]},y^{2,d}_{[0,4]},y^{3,d}_{[0,4]}\bigr).
		\end{align*}
		
\section{Conclusions}\label{sec5}
	We have given easily verifiable geometric conditions which assure that a selection of time derivatives of a flat output can be introduced as new (closed-loop) input, and have shown how to systematically construct the feedback required for actually introducing these derivatives as new input. Subsequently, we have proven in a finite-dimensional geometric framework that for every~${(x,u)}$-flat output \eqref{eq:xu_flat_output_procedure} of a system \eqref{eq:intro_nlsys} there exists an $m$-tuple $\kappa$ with $|\kappa|=n$ such that~${v=y_{[\kappa]}}$ is a feasible input, and that such an input can be introduced by a quasi-static feedback of the state $x$. Compared to the well-known exact linearization by a quasi-static feedback of a generalized Brunovsk\'y state, this approach has the advantage that it requires only knowledge of the state $x$ and not of time derivatives of the flat output up to a certain order. Furthermore, we have shown that on the basis of such an exact feedback linearization it is possible to systematically design a tracking control which only depends on the state $x$ and the reference trajectory, i.e., again without using a generalized Brunovsk\'y state. Future research will address the exact feedback linearization with respect to general flat outputs \eqref{eq:y}, which may also depend on time derivatives of the input $u$. The main challenge in generalizing the proposed results to flat outputs \eqref{eq:y} consist in finding a generalization or an alternative to the algorithm proposed in Section \ref{subsec:exact_lin_xu} for finding an $m$-tuple $\kappa$ with $|\kappa|=n$ such that $v=y_{[\kappa]}$ is a~feasible input. A system which is indeed flat but not $(x,u)$-flat can be found in \cite{GstottnerKolarSchoberl:2023}.

\subsection*{Acknowledgements}
This research was funded in whole, or in part, by the Austrian Science Fund (FWF) P32151 and P36473.
The authors would also like to thank J.~Rudolph and A.~Irscheid for interesting discussions concerning the exact linearization by quasi-static feedback of generalized states. Furthermore, we would like to thank the anonymous referees whose detailed reviews helped to improve the final version of the paper.

\pdfbookmark[1]{References}{ref}
\LastPageEnding

\end{document}

%% file: crane.pdf_tex
\begingroup%
  \makeatletter%
  \providecommand\color[2][]{%
    \errmessage{(Inkscape) Color is used for the text in Inkscape, but the package 'color.sty' is not loaded}%
    \renewcommand\color[2][]{}%
  }%
  \providecommand\transparent[1]{%
    \errmessage{(Inkscape) Transparency is used (non-zero) for the text in Inkscape, but the package 'transparent.sty' is not loaded}%
    \renewcommand\transparent[1]{}%
  }%
  \providecommand\rotatebox[2]{#2}%
  \newcommand*\fsize{\dimexpr\f@size pt\relax}%
  \newcommand*\lineheight[1]{\fontsize{\fsize}{#1\fsize}\selectfont}%
  \ifx\svgwidth\undefined%
    \setlength{\unitlength}{708.66141732bp}%
    \ifx\svgscale\undefined%
      \relax%
    \else%
      \setlength{\unitlength}{\unitlength * \real{\svgscale}}%
    \fi%
  \else%
    \setlength{\unitlength}{\svgwidth}%
  \fi%
  \global\let\svgwidth\undefined%
  \global\let\svgscale\undefined%
  \makeatother%
  \begin{picture}(1,0.58)%
    \lineheight{1}%
    \setlength\tabcolsep{0pt}%
    \put(0,0){\includegraphics[width=\unitlength,page=1]{crane.pdf}}%
    \put(0.51057407,0.49180408){\color[rgb]{0,0,0}\makebox(0,0)[lt]{\lineheight{0}\smash{\begin{tabular}[t]{l}$x$\end{tabular}}}}%
    \put(0.34891074,0.50623395){\color[rgb]{0,0,0}\makebox(0,0)[lt]{\lineheight{0}\smash{\begin{tabular}[t]{l}$y$\end{tabular}}}}%
    \put(0.44470807,0.44391113){\color[rgb]{0,0,0}\makebox(0,0)[lt]{\lineheight{0}\smash{\begin{tabular}[t]{l}$z$\end{tabular}}}}%
    \put(0,0){\includegraphics[width=\unitlength,page=2]{crane.pdf}}%
    \put(0.38441678,0.0387207){\color[rgb]{0,0,0}\makebox(0,0)[lt]{\lineheight{0}\smash{\begin{tabular}[t]{l}$m_L$\end{tabular}}}}%
    \put(0.3475563,0.19324146){\color[rgb]{0,0,0}\makebox(0,0)[lt]{\lineheight{0}\smash{\begin{tabular}[t]{l}$\alpha$\end{tabular}}}}%
    \put(0.33883105,0.14721069){\color[rgb]{0,0,0}\makebox(0,0)[lt]{\lineheight{0}\smash{\begin{tabular}[t]{l}$\beta$\end{tabular}}}}%
    \put(0.54198389,0.32554492){\color[rgb]{0,0,0}\makebox(0,0)[lt]{\lineheight{0}\smash{\begin{tabular}[t]{l}$f_x$\end{tabular}}}}%
    \put(0.42857913,0.2133891){\color[rgb]{0,0,0}\makebox(0,0)[lt]{\lineheight{0}\smash{\begin{tabular}[t]{l}$f_y$\end{tabular}}}}%
    \put(0.55015045,0.52550049){\color[rgb]{0,0,0}\makebox(0,0)[lt]{\lineheight{0}\smash{\begin{tabular}[t]{l}$x_T$\end{tabular}}}}%
    \put(0.268694,0.5158548){\color[rgb]{0,0,0}\makebox(0,0)[lt]{\lineheight{0}\smash{\begin{tabular}[t]{l}$y_T$\end{tabular}}}}%
    \put(0,0){\includegraphics[width=\unitlength,page=3]{crane.pdf}}%
    \put(0.58897083,0.37461523){\color[rgb]{0,0,0}\makebox(0,0)[lt]{\lineheight{0}\smash{\begin{tabular}[t]{l}$m_T$\end{tabular}}}}%
    \put(0.34520929,0.0879895){\color[rgb]{0,0,0}\makebox(0,0)[lt]{\lineheight{0}\smash{\begin{tabular}[t]{l}$l$\end{tabular}}}}%
    \put(0,0){\includegraphics[width=\unitlength,page=4]{crane.pdf}}%
    \put(0.72662848,0.27254449){\color[rgb]{0,0,0}\makebox(0,0)[lt]{\lineheight{0}\smash{\begin{tabular}[t]{l}$m_B$\end{tabular}}}}%
    \put(0,0){\includegraphics[width=\unitlength,page=5]{crane.pdf}}%
  \end{picture}%
\endgroup%